\documentclass[10pt]{article}

\usepackage{amsfonts}
\usepackage{amssymb}
\usepackage{amsmath}
\usepackage{amsthm}
\usepackage{pst-node}
\usepackage{pst-plot}
\usepackage{epsfig}
\usepackage{german}
\usepackage{psfrag}

\def\negthickspace{\!\!\!}
\def\beltrami{\mathrel{\mathop{\mbox{$\triangle$}}\limits_{{\parbox[t][0cm]{0.3cm}{\vspace*{-0.25cm}\rule{0.3cm}{0.025cm}}}}}}

\newtheorem{Lem}{Lemma}{\alph{enumi}}
\newenvironment{lemma}{\begin{Lem}\hspace{-0.2cm} }{\end{Lem}}

\newtheorem{Theo}{Theorem}{\alph{enumi}}
\newenvironment{theorem}{\begin{Theo}\hspace{-0.2cm}}{\end{Theo}}

\newtheorem{Prop}{Proposition}{\alph{enumi}}
\newenvironment{proposition}{\begin{Prop}\hspace{-0.2cm}}{\end{Prop}}

\newtheorem{Cor}{Corollary}{\alph{enumi}}
\newenvironment{corollary}{\begin{Cor}\hspace{-0.2cm}}{\end{Cor}}

\newtheorem{Def}{Definition}{\alph{enumi}}
\newenvironment{definition}{\begin{Def}\hspace{-0.2cm}\rm }{\end{Def}}

\newcommand{\nicefrac}[2]
{\leavevmode \kern.1em\raise.5ex\hbox{\the\scriptfont0 #1}
             \kern-.1em/\kern-.15em\lower.25ex
             \hbox{\the\scriptfont0 #2}}

\def\bef#1{\parbox[c][0.7cm][c]{0.3cm}\centering{#1}}

\begin{document}

\begin{center}
{\Large{\bf Remarks on Nitsche's functional:}}\\[0.2cm]
{\Large{\bf The rotationally symmetric case}}\\[0.8cm]
{\large{\sc Steffen Fr\"ohlich}}\\[1cm]
{\small\bf Abstract}\\[0.4cm]
\begin{minipage}[c][2cm][l]{12cm}
{\small We investigate existence and stability of rotationally symmetric critical immersions of variational problems of higher order which were considered  in \cite{Nitsche_02} and \cite{Nitsche_03}.}\end{minipage}
\end{center}
\subsection{Introduction}
Let $\Omega\subset\mathbb R^2$ be a bounded and two-fold connected domain. We consider two-dimensional immersions
\begin{equation}\label{1.1}
  X=X(u,v)=(x^1(u,v),x^2(u,v),x^3(u,v))\in C^{4+\alpha}(\Omega,\mathbb R^3)\cap C^0(\overline\Omega,\mathbb R^3),
  \quad\alpha\in(0,1),
\end{equation}
with the property
\begin{equation}\label{1.2}
  W:=|X_u\wedge X_v|>0\quad\mbox{in}\ \Omega
\end{equation}
for the surface area element $W=W(u,v).$ Here, the indices $u$ and $v$ are the partial derivatives w.r.t. to the variables $u$ resp. $v,$ while $\wedge$ means the usual vector product in $\mathbb R^3.$\\[0.1cm]
Let $\kappa_i=\kappa_i(u,v),$ $i=1,2,$ denote the principle curvatures of the surface. Then, by
\begin{equation}\label{1.3}
  H(u,v):=\frac{\kappa_1(u,v)+\kappa_2(u,v)}{2}\,,\quad
  K(u,v):=\kappa_1(u,v)\kappa_2(u,v)
\end{equation}
we introduce its mean curvature and Gaussian curvature. We are concerned with rotationally symmetric critical immersions w.r.t. the variational problem
\begin{equation}\label{1.4}
  {\mathcal E}[X]
  :=\int\hspace{-0.25cm}\int\limits_{\hspace{-0.3cm}\Omega}
    (\alpha+\beta H^2-\gamma K)W\,dudv
  \longrightarrow\mbox{extr!}
\end{equation}
with positive constants $\alpha,\beta,\gamma\in\mathbb R.$\\[0.1cm]
For $\alpha\not=0,$ $\beta,\gamma=0,$ the functional ${\mathcal E}[X]$ is proportional to the classical area functional. Critical points of the accessory variational problem are minimal surfaces, in this case the catenoid. In section 2 we investigate the stability of the catenoid solution w.r.t. small perturbations. The methods used here are presented in more general context in \cite{Nitsche_01}.\\[0.1cm]
The case $\alpha,\gamma=0$ and $\beta\not=0$ leads to Willmore's functional for which we present selected numerical results. Further, we refer the reader to the textbook \cite{Willmore_01}, in particular chapter 7.\\[0.1cm]
{\it Given two coaxial circular boundary curves $\Gamma_1,\Gamma_2\in\mathbb R^3$ with common radius $R>0$ of distance $d>0,$ we consider rotationally symmetric critical points of (\ref{1.4})} encouraged by the J.C.C. Nitsche's treatises \cite{Nitsche_02} and \cite{Nitsche_03}.
\subsection{The catenoid}
\setcounter{equation}{0}
\setcounter{Lem}{0}
\setcounter{Prop}{0}
\setcounter{Def}{0}
\setcounter{Theo}{0}
\setcounter{Cor}{0}
\subsubsection{The catenary curve}
Let $I:=[x_\ell,x_r]\subset\mathbb R,$ $|I|:=|x_r-x_\ell|>0.$ A critical point $f\in C^{2+\alpha}(I,\mathbb R),$ $\alpha\in(0,1),$ of the variational problem
\begin{equation}\label{2.1}
  {\mathcal A}[f]:=2\pi\int\limits_If(x)\sqrt{1+f'(x)^2}\,dx\longrightarrow\mbox{extr!}
\end{equation}
is a solution of the non-linear Euler-Lagrange differential equation
\begin{equation}\label{2.2}
  f(x)f''(x)=1+f'(x)^2\,,\quad x\in(x_\ell,x_r),
\end{equation}
with boundary conditions $y_\ell=f(x_\ell)$ and $y_r=f(x_r).$ Note that $f$ and $f''$ has no zeros. With suitable integration constants $c_1,c_2\in\mathbb R$ we get the catenary solution
\begin{equation}\label{2.3}
  f(x)=c_1\cosh\left(\frac{x}{c_1}+c_2\right),\quad x\in I.
\end{equation}
Detailed calculations can be found e.g. in \cite{Fomin_Gelfand_01}, section 4, and in \cite{Bliss_01}, chapter IV.
\subsubsection{Stability of graphs}
{\it (a) Perturbation of the catenary curve}\\[0.2cm]
We investigate the stability of the catenary curve $f$ w.r.t. small perturbations. Let $\psi\in C^{2+\alpha}(I,\mathbb R),$ $\alpha\in(0,1),$ be given such that $f+\psi$ solves also the minimal surface equation (\ref{2.2}), that is
\begin{equation}\label{2.4}
  (f+\psi)(f''+\psi'')-1-(f'+\psi')^2=0\quad\mbox{in}\ I.
\end{equation}
We arrive at the non-linear differential equation
\begin{equation}\label{2.5}
  \psi''(x)-\frac{2f'(x)}{f(x)}\,\psi'(x)+\frac{f''(x)}{f(x)}\,\psi(x)
  =\frac{1}{f}\,\psi'(x)^2-\frac{1}{f}\,\psi(x)\psi''(x).
\end{equation}
Due to $ff''=1+f'^2>0$ the linear differential operator on the right hand side does not obey the maximum principle (see Proposition 6.2 of the appendix). Thus, we assume the existence of a positive stability function $\chi\in C^{2+\alpha}(I,\mathbb R)$ such that
\begin{equation}\label{2.6}
  \chi''(x)-\frac{2f'(x)}{f(x)}\,\chi'(x)+\frac{f''(x)}{f(x)}\,\chi(x)\le 0\quad\mbox{in}\ I,\quad
  \chi>0\quad\mbox{in}\ I.
\end{equation}
With the product trick $\psi=\varphi\chi$ we calculate
\begin{equation}\label{2.7}
  \varphi''(x)
  +2\left(\frac{\chi'}{\chi}-\frac{f'}{f}\right)\varphi'(x)
  +\left(\frac{\chi''}{\chi}-\frac{2f'}{f}\,\frac{\chi'}{\chi}+\frac{f''}{f}\right)\varphi(x)
  =\Phi(\varphi;\chi)
\end{equation}
with the non-linear right hand side
\begin{equation}\label{2.8}
  \Phi(\varphi;\chi)
  :=\frac{\chi}{f}\,\varphi'^2-\frac{\chi}{f}\,\varphi\varphi''
    +\left(\frac{\chi'^2}{f\chi}-\frac{\chi''}{f}\right)\varphi^2\,,\quad
  \Phi(0;\chi)=0.
\end{equation}
Now, the left hand side differential operator in (\ref{2.7}) obeys the maximum principle. Corollary 6.3 ensures the uniqueness of the solution of the below boundary value problem (\ref{2.20}).\\[0.2cm]
{\it (b) Schauder norms}
\begin{definition}
We introduce the norms
\begin{equation}\label{2.9}
\begin{array}{lll}
  \|u\|_{0,I}\negthickspace
  & := & \negthickspace\displaystyle
         \max_{x\in I}|u(x)|,\quad
         \|u\|_{1,I}
         \,:=\,\|u\|_{0,I}+\max_{x\in I}|u'(x)|, \\[0.4cm]
         \|u\|_{2,I}\negthickspace
  & := & \negthickspace\displaystyle
         \|u\|_{1,I}+\max_{x\in I}|u''(x)|,
\end{array}
\end{equation}
furthermore the H\"older norms
\begin{equation}\label{2.10}
  \|u\|_{k+\alpha,I}
  :=\|u\|_{k,I}
    +\max_{\genfrac{}{}{0pt}{1}{x_1,x_2\in I}{x_1\not=x_2}}\frac{|u^{(k)}(x_1)-u^{(k)}(x_2)|}{|x_1-x_2|^\alpha}\,,
  \quad k=0,1,2,
\end{equation}
as well as the semi norms
\begin{equation}\label{2.11}
  [u]_{k,I}:=\max_{x\in I}|u^{(k)}(x)|,\quad
  [u]_{k+\alpha,I}
  :=\max_{\genfrac{}{}{0pt}{1}{x_1,x_2\in I}{x_1\not=x_2}}\frac{|u^{(k)}(x_1)-u^{(k)}(x_2)|}{|x_1-x_2|^\alpha}
  \quad\mbox{for}\ k=0,1,2.
\end{equation}
\end{definition}
\noindent
{\sl (c) Contraction property of $\Phi(\eta;\chi)$}
\begin{lemma}
It holds
\begin{equation}\label{2.12}
  \|\Phi(\eta_1;\chi)-\Phi(\eta_2;\chi)\|_{\alpha,I}
  \le {\mathcal C}_1(\|\eta_1\|_{2+\alpha,I}+\|\eta_2\|_{2+\alpha,I})\|\eta_1-\eta_2\|_{2+\alpha,I}
\end{equation}
for all $\eta_1,\eta_2\in C^{2+\alpha}(I,\mathbb R)$ with the constant
\begin{equation}\label{2.13}
  {\mathcal C}_1
  :=\Big\{
      2\|\chi f^{-1}\|_{\alpha,I}
      +\|\chi'^2\chi^{-1}f^{-1}-\chi''f^{-1}\|_{\alpha,I}(1+|I|^{1-\alpha})^2
    \Big\}
    (1+|I|^{1-\alpha})^2\,.
\end{equation}
\end{lemma}
\begin{proof}
We calculate
\begin{equation}\label{2.14}
\begin{array}{lll}
  \Phi(\eta_1;\chi)-\Phi(\eta_2;\chi)\negthickspace
  & = & \negthickspace\displaystyle
        \frac{\chi}{f}\,(\eta_1'^2-\eta_2'^2)
        -\frac{\chi}{f}\,(\eta_1\eta_1''-\eta_2\eta_2'')
        -\left(\frac{\chi'^2}{f\chi}-\frac{\chi''}{f}\right)(\eta_1^2-\eta_2^2) \\[0.6cm]
  & = & \negthickspace\displaystyle
        \frac{\chi}{f}\,(\eta_1'+\eta_2')(\eta_1'-\eta_2')
        -\frac{\chi}{f}\,\eta_1(\eta_1''-\eta_2'')
        -\frac{\chi}{f}\,\eta_2''(\eta_1-\eta_2) \\[0.6cm]
  &   & \negthickspace\displaystyle
        -\left(\frac{\chi'^2}{f\chi}-\frac{\chi''}{f}\right)(\eta_1+\eta_2)(\eta_1-\eta_2).
\end{array}
\end{equation}
Now, we estimate as follows:
\begin{equation}\label{2.15}
\begin{array}{l}
  \|\Phi(\eta_1;\chi)-\Phi(\eta_2;\chi)\|_{\alpha,I} \\[0.2cm]
  \hspace*{0.6cm}\displaystyle
  \le\,\|\chi f^{-1}\|_{\alpha,I}\|\eta'_1+\eta'_2\|_{\alpha,I}\|\eta'_1-\eta'_2\|_{\alpha,I} \\[0.3cm]
  \hspace*{1.2cm}\displaystyle
    +\,\|\chi f^{-1}\|_{\alpha,I}\|\eta_1\|_{\alpha,I}\|\eta''_1-\eta''_2\|_{\alpha,I}
    +\|\chi f^{-1}\|_{\alpha,I}\|\eta''_2\|_{\alpha,I}\|\eta_1-\eta_2\|_{\alpha,I} \\[0.3cm]
  \hspace*{1.2cm}\displaystyle
    +\,\|\chi'^2\chi^{-1}f^{-1}-\chi''f^{-1}\|_{\alpha,I}
       \|\eta_1+\eta_2\|_{\alpha,I}\|\eta_1-\eta_2\|_{\alpha,I} \\[0.3cm]
  \hspace*{0.6cm}\displaystyle
  \le\,\|\chi f^{-1}\|_{\alpha,I}(\|\eta_1\|_{1+\alpha,I}+\|\eta_2\|_{1+\alpha,I})\|\eta_1-\eta_2\|_{1+\alpha,I} \\[0.3cm]
  \hspace*{1.2cm}\displaystyle
    +\,\|\chi f^{-1}\|_{\alpha,I}\|\eta_1\|_{\alpha,I}\|\eta_1-\eta_2\|_{2+\alpha,I}
    +\|\chi f^{-1}\|_{\alpha,I}\|\eta_2\|_{2+\alpha,I}\|\eta_1-\eta_2\|_{\alpha,I} \\[0.3cm]
  \hspace*{1.2cm}\displaystyle
    +\,\|\chi'^2\chi^{-1}f^{-1}-\chi''f^{-1}\|_{\alpha,I}(\|\eta_1\|_{\alpha,I}+\|\eta_2\|_{\alpha,I})
     \|\eta_1-\eta_2\|_{\alpha,I}\,.
\end{array}
\end{equation}
Note that
\begin{equation}\label{2.16}
\begin{array}{lll}
  \|\eta\|_{\alpha,I}\negthickspace
  & \le & \negthickspace\displaystyle
          \|\eta\|_{0,I}
          +\max_{\genfrac{}{}{0pt}{1}{x_1,x_2\in I}{x_1\not=x_2}}
           \frac{|\eta(x_1)-\eta(x_2)|}{|x_1-x_2|}\,|I|^{1-\alpha} \\[0.7cm]
  & \le & \negthickspace\displaystyle
          \|\eta\|_{0,I}
          +|I|^{1-\alpha}\max_{\genfrac{}{}{0pt}{1}{x_1,x_2\in I}{x_1\not=x_2}}
           \max_{\xi\in[x_1,x_2]}|\eta'(\xi)| \\[0.7cm]
  & \le & \negthickspace\displaystyle
          \|\eta\|_{0,I}
          +|I|^{1-\alpha}\max_{\genfrac{}{}{0pt}{1}{x_1,x_2\in I}{x_1\not=x_2}}\|\eta\|_{1,[x_1,x_2]} \\[0.7cm]
  & \le & \negthickspace\displaystyle
          (1+|I|^{1-\alpha})\|\eta\|_{1+\alpha,I}\,.
\end{array}
\end{equation}
In the same way we have
\begin{equation}\label{2.17}
  \|\eta\|_{1+\alpha,I}
  \le(1+|I|^{1-\alpha})\|\eta\|_{2+\alpha,I}\,,\quad
  \|\eta\|_{\alpha,I}\le(1+|I|^{1-\alpha})^2\|\eta\|_{2+\alpha,I}\,.
\end{equation}
We insert these inequalities into (\ref{2.15}) to get
\begin{equation}\label{2.18}
\begin{array}{l}
  \|\Phi(\eta_1;\chi)-\Phi(\eta_2;\chi)\|_{\alpha,I}\\[0.3cm]
  \hspace*{0.6cm}\displaystyle
  \le\,\|\chi f^{-1}\|_{\alpha,I}(1+|I|^{1-\alpha})^2
       (\|\eta_1\|_{2+\alpha,I}+\|\eta_2\|_{2+\alpha,I})\|\eta_1-\eta_2\|_{2+\alpha,I} \\[0.3cm]
  \hspace*{1.2cm}\displaystyle
    +\,\|\chi f^{-1}\|_{\alpha,I}(1+|I|^{1-\alpha})^2\|\eta_1\|_{2+\alpha,I}\|\eta_1-\eta_2\|_{2+\alpha,I} \\[0.3cm]
  \hspace*{1.2cm}\displaystyle
    +\,\|\chi f^{-1}\|_{\alpha,I}(1+|I|^{1-\alpha})^2\|\eta_2\|_{2+\alpha,I}\|\eta_1-\eta_2\|_{2+\alpha,I} \\[0.3cm]
  \hspace*{1.2cm}\displaystyle
    +\,\|\chi'^2\chi^{-1}f^{-1}-\chi''f^{-1}\|_{\alpha,I}(1+|I|^{1-\alpha})^4
       (\|\eta_1\|_{2+\alpha,I}+\|\eta_2\|_{2+\alpha,I})\|\eta_1-\eta_2\|_{2+\alpha,I}\,.
\end{array}
\end{equation}
This is the statement.
\end{proof}
\noindent
For $\varphi\in C^{2+\alpha}(I,\mathbb R)$ we define the linear differential operator
\begin{equation}\label{2.19}
  {\mathcal L}[\varphi]
  :=\varphi''(x)
    +2\left(\frac{\chi'}{\chi}-\frac{f'}{f}\right)\varphi'(x)
    +\left(\frac{\chi''}{\chi}-\frac{2f'\chi'}{f\chi}+\frac{f''}{f}\right)\varphi(x).
\end{equation}
Due to (\ref{2.6}) it obeys the maximum principle.\\[0.1cm]
Successively we will solve the boundary value problem (let $\Phi(\varphi):=\Phi(\varphi;\chi)$)
\begin{equation}\label{2.20}
\begin{array}{l}
  {\mathcal L}[\varphi]=\Phi(\varphi)\quad\mbox{in}\ I, \\[0.2cm]
  \varphi(x_\ell)=\varphi_\ell\,,\ \varphi(x_r)=\varphi_r
\end{array}
\end{equation}
with the non-linear right hand side (\ref{2.8}) and given boundary data $\varphi(x_\ell)=\varphi_\ell$ and $\varphi(x_r)=\varphi_r.$ For this we start with the function $\varphi_0\equiv 0$ and consider the linear problems
\begin{equation}\label{2.21}
\begin{array}{l}
  {\mathcal L}[\varphi_k]=\Phi(\varphi_{k-1})\quad\mbox{in}\ I, \\[0.2cm]
  \varphi(x_\ell)=\varphi_\ell\,,\ \varphi(x_r)=\varphi_r\quad\mbox{for}\ k=1,2,\ldots
\end{array}
\end{equation}
From ${\mathcal L}[\varphi_1]=\Phi(0)=0$ and the Schauder estimate
\begin{equation}\label{2.22}
  \|\varphi_k\|_{2+\alpha,I}
  \le{\mathcal C}_2\|\Phi(\varphi_{k-1})\|_{\alpha,I}
     +{\mathcal C}_3\max\{|\varphi(x_\ell)|,|\varphi(x_r)|\}
\end{equation}
(see the global $C^{2+\alpha}$-estimate (\ref{6.40}) together with the $C^0$-estimate (\ref{6.7}) and the a priori constants ${\mathcal C}_2,{\mathcal C}_3\in(0,+\infty)$ following from it) we conclude
\begin{equation}\label{2.23}
  \|\varphi_1\|_{2+\alpha,I}
  \le{\mathcal C}_3\max\{|\varphi(x_\ell)|,|\varphi(x_r)|\}
  =:{\mathcal C}_3a,
\end{equation}
where $a:=\max\{|\varphi(x_\ell)|,|\varphi(x_r)|\}.$ Using (\ref{2.12}) we get
\begin{equation}\label{2.24}
\begin{array}{lll}
  \|\varphi_{k+1}-\varphi_k\|_{2+\alpha,I}\negthickspace
  & \le & \negthickspace\displaystyle
          {\mathcal C}_2\|\Phi(\varphi_k)-\Phi(\varphi_{k-1})\|_{\alpha,I} \\[0.2cm]
  & \le & \negthickspace\displaystyle
          {\mathcal C}_1{\mathcal C}_2(\|\varphi_k\|_{2+\alpha,I}+\|\varphi_{k-1}\|_{2+\alpha,I})
          \|\varphi_k-\varphi_{k-1}\|_{2+\alpha,I}\,.
\end{array}
\end{equation}
For given $\varepsilon\in(0,1)$ we choose $a$ sufficiently small such that
\begin{equation}\label{2.25}
  2a\,{\mathcal C}_1{\mathcal C}_2{\mathcal C}_3
  \left(1+\frac{a\,{\mathcal C}_1{\mathcal C}_2{\mathcal C}_3}{1-\varepsilon}\right)\le\varepsilon.
\end{equation}
Then, by induction we prove
\begin{equation}\label{2.26}
  \|\varphi_n\|_{2+\alpha,I}
  \le a\,{\mathcal C}_3\left(1+\frac{a\,{\mathcal C}_1{\mathcal C}_2{\mathcal C}_3}{1-\varepsilon}\right)
  \quad\mbox{for}\ n=1,2,\ldots
\end{equation}
First, it holds $\|\varphi_1\|_{2+\alpha,I}\le a{\mathcal C}_3.$ Furthermore, (\ref{2.24}) gives
\begin{equation}\label{2.27}
\begin{array}{lll}
  \|\varphi_2\|_{2+\alpha,I}\negthickspace
  & \le & \negthickspace\displaystyle
          \|\varphi_2-\varphi_1\|_{2+\alpha,I}+\|\varphi_1\|_{2+\alpha,I}
          \,\le\,{\mathcal C}_1{\mathcal C}_2\|\varphi_1\|_{2+\alpha,I}^2+\|\varphi_1\|_{2+\alpha,I} \\[0.2cm]
  & \le & \negthickspace\displaystyle
          a\,{\mathcal C}_3(1+a\,{\mathcal C}_1{\mathcal C}_2{\mathcal C}_3).
\end{array}
\end{equation}
Thus, the statement is true for $n=1$ und $n=2.$ Let it be proved for $n=1,2,\ldots,m.$ For $k\le m$ we calculate (see (\ref{2.24}))
\begin{equation}\label{2.28}
  \|\varphi_{k+1}-\varphi_k\|_{2+\alpha,I}
  \le 2a\,{\mathcal C}_1{\mathcal C}_2{\mathcal C}_3
      \left(1+\frac{a\,{\mathcal C}_1{\mathcal C}_2{\mathcal C}_3}{1-\varepsilon}\right)
      \|\varphi_k-\varphi_{k-1}\|_{2+\alpha,I}
  \le \varepsilon\|\varphi_k-\varphi_{k-1}\|_{2+\alpha,I}\,.
\end{equation}
It follows that (see (\ref{2.27}))
\begin{equation}\label{2.29}
\begin{array}{lll}
  \|\varphi_{m+1}\|_{2+\alpha,I}\negthickspace
  & \le & \negthickspace\displaystyle
          \|\varphi_1\|_{2+\alpha,I}
          +\sum_{k=1}^m\|\varphi_{k+1}-\varphi_k\|_{2+\alpha,I} \\[0.6cm]
  & \le & \negthickspace\displaystyle
          \|\varphi_1\|_{2+\alpha,I}
          +(\varepsilon^{m-1}+\varepsilon^{m-2}+\ldots+\varepsilon+1)\|\varphi_2-\varphi_1\|_{2+\alpha,I} \\[0.4cm]
  & \le & \negthickspace\displaystyle
          a\,{\mathcal C}_3+\frac{1}{1-\varepsilon}\,a^2{\mathcal C}_1{\mathcal C}_2{\mathcal C}_3^2
          \,=\,a\,{\mathcal C}_3\left(1+\frac{a\,{\mathcal C}_1{\mathcal C}_2{\mathcal C}_3}{1-\varepsilon}\right).
\end{array}
\end{equation}
This proves (\ref{2.26}) and we arrive at
\begin{equation}\label{2.30}
  \|\varphi_{k+1}-\varphi_k\|_{2+\alpha,I}\le\varepsilon\|\varphi_k-\varphi_{k-1}\|_{2+\alpha,I}
  \quad\mbox{for all}\ k=1,2,\ldots
\end{equation}
for $\varepsilon\in(0,1).$ By Banach's fix point theorem the sequence $\{\varphi_k\}_{k=0,1,\ldots}$ converges uniformly in $C^{2+\alpha}$ to a function $\varphi\in C^{2+\alpha}(I,\mathbb R)$ which solves the boundary value problem (\ref{2.20}). Via $\psi=\varphi\chi$ we have found also a solution of (\ref{2.5}) with sufficiently small boundary data.
\begin{theorem}
Let $f\in C^{2+\alpha}(I,\mathbb R)$ be a solution of the minimal surface equation (\ref{2.2}) and let it be stable in the sense of (\ref{2.6}). For sufficiently small $\varepsilon>0$ there exists $\psi\in C^{2+\alpha}(I,\mathbb R)$ with boundary values $|\psi(x_\ell)|\le\varepsilon$ and $|\psi(x_r)|\le\varepsilon$ such that $f+\psi$ solves also (\ref{2.2}).
\end{theorem}
\subsection{Variational problems of higher order}
\subsubsection{The first variation}
\setcounter{equation}{0}
\setcounter{Lem}{0}
\setcounter{Prop}{0}
\setcounter{Def}{0}
\setcounter{Theo}{0}
\setcounter{Cor}{0}
Now we address the first variation of the functional ${\mathcal E}[X]$ from (\ref{1.4}). We introduce conformal parameters $(u,v)\in\Omega$ with the properties
\begin{equation}\label{3.1}
  |X_u|^2=W=|X_v|^2\,,\quad
  X_u\cdot X_v^t=0
  \quad\mbox{in}\ \Omega.
\end{equation}
We consider the perturbation
\begin{equation}\label{3.2}
  \widetilde X(u,v):=X(u,v)+\varepsilon\{\varphi(u,v)X_u(u,v)+\psi(u,v)X_v(u,v)+\chi(u,v)N(u,v)\},
  \quad(u,v)\in\overline\Omega,
\end{equation}
where $\varepsilon\in(-\varepsilon_0,+\varepsilon_0).$ Let
\begin{equation}\label{3.3}
  Z(u,v):=\varphi(u,v)X_u(u,v)+\psi(u,v)X_v(u,v)+\chi(u,v)N(u,v).
\end{equation}
From \cite{Nitsche_03} and \cite{Mosel_01} we have
\begin{lemma}
For the first variation w.r.t. the vector field $Z=Z(u,v)$ there hold
\begin{equation}\label{3.4}
\begin{array}{rcl}
  \displaystyle
  \delta\int\hspace{-0.25cm}\int\limits_{\hspace{-0.3cm}\Omega}W\,dudv\negthickspace
  & = & \displaystyle\negthickspace
        -2\int\hspace{-0.25cm}\int\limits_{\hspace{-0.3cm}\Omega}H(Z\cdot N^t)W\,dudv
        -\int\limits_{\partial\Omega}[Z,N,X']\,ds, \\[0.8cm]
  \displaystyle
  \delta\int\hspace{-0.25cm}\int\limits_{\hspace{-0.3cm}\Omega}KW\,dudv\negthickspace
  & = & \displaystyle\negthickspace
        -\int\limits_{\partial\Omega}K[Z,N,X']\,ds
        +\int\limits_{\partial\Omega}
           \left\{
             \Lambda^{(1)}\frac{\partial(Z\cdot N^t)}{\partial\nu}-\Lambda^{(2)}(Z\cdot N^t)
             \right\}ds, \\[0.8cm]
  \displaystyle
  \delta\int\hspace{-0.25cm}\int\limits_{\hspace{-0.3cm}\Omega}H^2W\,dudv\negthickspace
  & = & \displaystyle\negthickspace
        \int\hspace{-0.25cm}\int\limits_{\hspace{-0.3cm}\Omega}
        \{\triangle H+2H(H^2-K)W\}(Z\cdot N^t)\,dudv \\[0.8cm]
  &   & \displaystyle\negthickspace
        -\,\int\limits_{\partial\Omega}H^2[Z,N,X']\,ds
        +\int\limits_{\partial\Omega}
         \left\{
           H\frac{\partial(Z\cdot N^t)}{\partial\nu}-(Z\cdot N^t)\frac{\partial H}{\partial\nu}
         \right\}ds
\end{array}
\end{equation}
with the settings
\begin{equation}\label{3.5}
\begin{array}{lll}
  \Lambda^{(1)}\negthickspace
  & = & \negthickspace
        \kappa_n\,,\\[0.2cm]
  \Lambda^{(2)}\negthickspace
  & = & \negthickspace\displaystyle
        \frac{\partial\tau}{\partial s}
        +\frac{\partial}{\partial s}\,\bigg(\frac{\kappa_n}{\kappa^2}\,\frac{\partial\kappa_g}{\partial s}\bigg)
        -\frac{\partial}{\partial s}\,\bigg(\frac{\kappa_g}{\kappa^2}\,\frac{\partial\kappa_n}{\partial s}\bigg).
\end{array}
\end{equation}
Here, $\kappa$ and $\tau$ are the curvature and the torsion of the boundary curves, while $\kappa_n$ and $\kappa_g$ mean their normal curvature and geodesic curvature, resp., w.r.t. the surface $X=X(u,v).$
\end{lemma}
\noindent
The Euler-Lagrange differential equation
\begin{equation}\label{3.6}
  \beta\{\triangle H+2H(H^2-K)W\}-2\alpha HW=0
\end{equation}
has to be coupled with the natural boundary conditions
\begin{equation}\label{3.7}
\begin{array}{lll}
  0 \negthickspace
  & = & \displaystyle\negthickspace
        \int\limits_{\partial\Omega}\{\gamma K-\beta H^2-\alpha\}[Z,N,X']\,ds
        +\int\limits_{\partial\Omega}
         \left\{
           \gamma\Lambda^{(2)}-\beta\frac{\partial H}{\partial\nu}
         \right\}(Z\cdot N^t)\,ds \\[0.8cm]
  &   & \displaystyle\negthickspace
        +\,\int\limits_{\partial\Omega}
           \{\beta H-\gamma\Lambda^{(1)}\}
           \frac{\partial(Z\cdot N^t)}{\partial\nu}\,ds.
\end{array}
\end{equation}
Furthermore, we add the mean curvature differential system
\begin{equation}\label{3.8}
  \triangle X=2H(X)X_u\wedge X_v\,.
\end{equation}
In the case of fixed boundary conditions there hold $\chi|_{\partial\Omega}=0$ and $(Z\wedge X')|_{\partial\Omega}=0,$ that is
\begin{equation}\label{3.9}
  \beta H-\gamma\Lambda^{(1)}=0
\end{equation}
(for a detailed discussion we refer to \cite{Nitsche_02} and \cite{Nitsche_01}.)
\subsubsection{Differential equations for graph solutions}
Let $\beltrami$ denote the Laplace-Beltrami operator w.r.t. the metric $ds^2=|X_u|^2\,du^2+2X_u\cdot X_v^t\,dudv+|X_v|^2\,dv^2.$ If $\beta\not=0$ we can rewrite (\ref{3.6}) and (\ref{3.8}) to get
\begin{equation}\label{3.10}
  \beltrami\hspace{-0.08cm}H+2H(H^2-K)-\frac{2\alpha}{\beta}\,H=0,\quad
  W\hspace{-0.1cm}\beltrami\hspace{-0.08cm}X=2H(X_u\wedge X_v).
\end{equation}
Now, assume the surface represents a graph $(x,y,z(x,y)).$ We calculate
\begin{equation}\label{3.11}
\begin{array}{lll}
  \displaystyle
  W\hspace{-0.1cm}\beltrami\hspace{-0.04cm}z\negthickspace
  & = & \negthickspace\displaystyle
        \frac{\partial}{\partial x}\bigg(\frac{1+z_y^2}{W}\,z_x\bigg)
        -\frac{\partial}{\partial x}\bigg(\frac{z_xz_y}{W}\,z_y\bigg)
        -\frac{\partial}{\partial y}\bigg(\frac{z_xz_y}{W}\,z_x\bigg)
        +\frac{\partial}{\partial y}\bigg(\frac{1+z_x^2}{W}\,z_y\bigg) \\[0.6cm]
  & = & \negthickspace\displaystyle
        \frac{1}{W^3}\,\Big\{(1+z_y^2)z_{xx}-2z_xz_yz_{xy}+(1+z_x^2)z_{yy}\Big\}.
\end{array}
\end{equation}
Thus, it holds
\begin{equation}\label{3.12}
  (1+z_y^2)z_{xx}-2z_xz_yz_{xy}+(1+z_x^2)z_{yy}=2H(1+z_x^2+z_y^2)^\frac{3}{2}\,.
\end{equation}
Analogously we have
\begin{equation}\label{3.13}
\begin{array}{lll}
  W\hspace*{-0.06cm}\beltrami\hspace*{-0.1cm}H\negthickspace
  & = & \negthickspace\displaystyle
        \frac{2z_yz_{xy}}{W}\,H_x
        -\frac{(1+z_y^2)(z_xz_{xx}+z_yz_{xy})}{W^3}\,H_x
        +\frac{1+z_y^2}{W}\,H_{xx} \\[0.6cm]
  &   & \negthickspace\displaystyle
        -\,\frac{z_yz_{xx}+z_xz_{xy}}{W}\,H_y
        +\frac{z_xz_y(z_xz_{xx}+z_yz_{xy})}{W^3}\,H_y
        -\frac{z_xz_y}{W}\,H_{xy} \\[0.6cm]
  &   & \negthickspace\displaystyle
        -\,\frac{z_yz_{xy}+z_xz_{yy}}{W}\,H_x
        +\frac{z_xz_y(z_xz_{xy}+z_yz_{yy})}{W^3}\,H_x
        -\frac{z_xz_y}{W}\,H_{xy} \\[0.6cm]
  &   & \negthickspace\displaystyle
        +\,\frac{2z_xz_{xy}}{W}\,H_y
        -\frac{(1+z_x^2)(z_xz_{xy}+z_yz_{yy})}{W^3}\,H_y
        +\frac{1+z_x^2}{W}\,H_{yy}
\end{array}
\end{equation}
for the mean curvature.
\subsubsection{The rotationally symmetric case}
We consider rotationally symmetric surfaces: The meridian function $f(x):=z(x,0),$ where $f(x)>0$ for all $x\in[x_\ell,x_r],$ rotates about the $x$-axis. We have
\begin{equation}\label{3.14}
\begin{array}{l}
  z_x(x,0)=f'(x),\quad z_y(x,0)=0, \\[0.2cm]
  \displaystyle
  z_{xx}(x,0)=f''(x),\quad z_{xy}(x,0)=0,\quad z_{yy}(x,0)=-\frac{1}{f(x)}
\end{array}
\end{equation}
as well as
\begin{equation}\label{3.15}
\begin{array}{l}
  H_x(x,0)=H'(x),\quad H_y(x,0)=0, \\[0.2cm]
  H_{xx}(x,0)=H''(x),\quad H_{xy}(x,0)=0,\quad H_{yy}(x,0)=0.
\end{array}
\end{equation}
The mean curvature equation reads as
\begin{equation}\label{3.16}
  f''=2H(1+f'^2)^\frac{3}{2}+\frac{1+f'^2}{f}
\end{equation}
while from (\ref{3.13}) we conclude
\begin{equation}\label{3.17}
  \beltrami\hspace*{-0.1cm}H
  =\frac{f'(x)}{1+f'(x)^2}\left(\frac{1}{f(x)}-\frac{f''(x)}{1+f'(x)^2}\right)H'(x)
   +\frac{1}{1+f'(x)^2}\,H''(x).
\end{equation}
Investing
\begin{equation}\label{3.18}
  K(x)=-\frac{f''(x)}{f(x)\{1+f'(x)^2\}}
\end{equation}
along the meridian curve, from (\ref{3.10}) and (\ref{3.17}) we get
\begin{equation}\label{3.19}
  H''+f'\left(\frac{1}{f}-\frac{f''}{1+f'^2}\right)H'
  +2H^3(1+f'^2)+2\,\frac{f''}{f}\,H-\frac{2\alpha}{\beta}\,(1+f'^2)H=0.
\end{equation}
\subsection{Successive approximation}
\setcounter{equation}{0}
\setcounter{Lem}{0}
\setcounter{Prop}{0}
\setcounter{Def}{0}
\setcounter{Theo}{0}
\setcounter{Cor}{0}
\subsubsection{Perturbation of the mean curvature equation}
We start with
\begin{equation}\label{4.1}
  H=\frac{f''}{2(1+f'^2)^\frac{3}{2}}-\frac{1}{2f\sqrt{1+f'^2}}\,.
\end{equation}
Let $f\in C^{2+\alpha}(I,\mathbb R)$ solve the minimal surface equation. For a perturbation $\psi\in C^{2+\alpha}(I,\mathbb R)$ we consider the mean curvature
\begin{equation}\label{4.2}
  \widetilde H
  =\frac{f''+\psi''}{2(1+(f'+\psi')^2)^\frac{3}{2}}
   -\frac{1}{2(f+\psi)\sqrt{1+(f'+\psi')^2}}
\end{equation}
along the meridian curve of the rotationally symmetric surface $f+\psi.$ Using
\begin{equation}\label{4.3}
\begin{array}{l}
  \displaystyle
  \frac{1}{f+\psi}
  \,=\,\frac{1}{f}-\frac{\psi}{f^2}+\ldots, \\[0.6cm]
  \displaystyle
  \frac{1}{\sqrt{1+f'^2+2f'\psi'+\psi'^2}}
  \,=\,\frac{1}{\sqrt{1+f'^2}}-\frac{f'\psi'}{(1+f'^2)^\frac{3}{2}}+\ldots, \\[0.6cm]
  \displaystyle
  \frac{1}{(1+f'^2+2f'\psi'+\psi'^2)^\frac{3}{2}}
  \,=\,\frac{1}{(1+f'^2)^\frac{3}{2}}-\frac{3f'\psi'}{(1+f'^2)^\frac{5}{2}}+\ldots
\end{array}
\end{equation}
we calculate
\begin{equation}\label{4.4}
\begin{array}{lll}
  \widetilde H\negthickspace
  & = & \negthickspace\displaystyle
        \frac{1}{2}
        \left[
          \frac{1}{(1+f'^2)^\frac{3}{2}}-\frac{3f'\psi'}{(1+f'^2)^\frac{5}{2}}
        \right](f''+\psi'') \\[0.8cm]
  &   & \negthickspace\displaystyle
        -\,\frac{1}{2}
           \left[
             \frac{1}{f}-\frac{\psi}{f^2}
           \right]
           \left[
             \frac{1}{\sqrt{1+f'^2}}-\frac{f'\psi'}{(1+f'^2)^\frac{3}{2}}
           \right]+\ldots \\[0.8cm]
  & = & \negthickspace\displaystyle
        \frac{1}{2\sqrt{1+f'^2}}
        \left[
          \frac{f''}{1+f'^2}-\frac{1}{f}
        \right]
        +\frac{1}{2(1+f'^2)^\frac{3}{2}}\,\psi''
        +\frac{f'}{2(1+f'^2)^\frac{3}{2}}
         \left[
           \frac{1}{f}-\frac{3f''}{1+f'^2}
         \right]\psi' \\[0.8cm]
  &   & \negthickspace\displaystyle
        +\,\frac{1}{2f^2\sqrt{1+f'^2}}\,\psi
        +\ldots
\end{array}
\end{equation}
for $\|\psi\|_{1,I}\le\varepsilon$ sufficiently small.
\goodbreak\noindent
Because $f$ solves the minimal surface equation (\ref{2.2}), we get
\begin{equation}\label{4.5}
  \widetilde H
  =\frac{1}{2(1+f'^2)^\frac{3}{2}}\,\psi''
   -\frac{f'}{f(1+f'^2)^\frac{3}{2}}\,\psi'
   +\frac{f''}{2f(1+f'^2)^\frac{3}{2}}\,\psi
   +\ldots
\end{equation}
\begin{lemma}
Let the perturbation $\psi\in C^{2+\alpha}(I,\mathbb R)$ be sufficiently small w.r.t. the $C^1$-norm. Let $f\in C^{2+\alpha}(I,\mathbb R)$ solve the minimal surface equation. Then it holds
\begin{equation}\label{4.6}
  \psi''-\frac{2f'}{f}\,\psi'+\frac{f''}{f}\,\psi
  =2(1+f'^2)^\frac{3}{2}\,\widetilde H+\Psi_1(\psi)
\end{equation}
with the mean curvature $\widetilde H$ along $f+\psi$ and the non-linear term $\Psi_1(\psi)$ which collects all super-linear terms for $\psi$ and its first and second derivatives. Moreover, with a stability function $\chi$ from (\ref{2.6}), the product $\varphi=\psi\chi^{-1}$ satisfies
\begin{equation}\label{4.7}
  \varphi''
  +2\left(
     \frac{\chi'}{\chi}-\frac{f'}{f}
   \right)\varphi'
  +\left(
     \frac{\chi''}{\chi}-\frac{2f'\chi'}{f\chi}+\frac{f''}{f}
   \right)\varphi
  =\Phi_1(\varphi,\widetilde H)
\end{equation}
with the non-linear right hand side
\begin{equation}\label{4.8}
  \Phi_1(\varphi,\widetilde H):=\frac{2(1+f'^2)^\frac{3}{2}}{\chi}\,\widetilde H+\chi^{-1}\Psi_1(\varphi\chi).
\end{equation}
\end{lemma}
\noindent
Let us denote the linear differential operator on the left hand side of (\ref{4.7}) with ${\mathcal L}_1.$ Then we consider the boundary value problem
\begin{equation}\label{4.9}
\begin{array}{l}
  {\mathcal L}_1[\varphi]=\Phi_1(\varphi,\widetilde H)\quad\mbox{in}\ I, \\[0.2cm]
  \varphi(x_\ell)=0,\ \varphi(x_r)=0.
\end{array}
\end{equation}
We supplement the boundary value problem (\ref{4.13}) for $\widetilde H.$
\subsubsection{Perturbation of the Euler-Lagrange equation}
Let the variation $f+\varphi,$ where $\|\varphi\|_{1,I}$ is sufficiently small and $f$ solves the minimal surface equation, satisfy (\ref{3.19}), that is
\begin{equation}\label{4.10}
\begin{array}{lll}
  0\negthickspace
  & = & \negthickspace\displaystyle
        \widetilde H''
        +(f'+\varphi')
        \left(
          \frac{1}{f+\varphi}
          -\frac{f''+\varphi''}{1+f'^2+2f'\varphi'+\varphi'^2}
        \right)
        \widetilde H'
        +2(1+f'^2+2f'\varphi'+\varphi'^2)\widetilde H^3 \\[0.8cm]
  &   & \negthickspace\displaystyle
        +\,2\,\frac{f''+\varphi''}{f+\varphi}\,\widetilde H
        -\frac{2\alpha}{\beta}\,(1+f'^2+2f'\varphi'+\varphi'^2)\widetilde H\,.
\end{array}
\end{equation}
\begin{lemma}
For the mean curvature along $f+\varphi$ it holds
\begin{equation}\label{4.11}
  \widetilde H''
  +\frac{2f''}{f}\,\widetilde H
  -\frac{2\alpha}{\beta}\,(1+f'^2)\widetilde H
  =\Phi_2(\varphi,\widetilde H)
\end{equation}
with the non-linear right hand side
\begin{equation}\label{4.12}
\begin{array}{lll}
  \Phi_2(\varphi,\widetilde H)\negthickspace
  & = & \negthickspace\displaystyle
        -2(1+f'^2)\widetilde H^3
        +(f'+\varphi')\left\{
           \frac{\varphi}{f^2}
           -\frac{2f'f''\varphi'}{(1+f'^2)^2}
           +\frac{\varphi''}{1+f'^2}
           -\frac{2f'\varphi'\varphi''}{(1+f'^2)^2}
         \right\}\widetilde H' \\[0.6cm]
  &   & \negthickspace\displaystyle
        -\,2(f'+\varphi')\varphi'\widetilde H^3
        +\left\{
           \frac{2f''\varphi}{f^2}
           -\frac{2\varphi''}{f}
           +\frac{2\varphi\varphi''}{f^2}
           +\frac{2\alpha}{\beta}\,(2f'\varphi'+\varphi'^2)
         \right\}\widetilde H
        +\ldots,
\end{array}
\end{equation}
where $\ldots$ means terms of higher order of $\psi$ and its derivatives.
\end{lemma}
\noindent
Let ${\mathcal L}_2$ denote the linear differential operator on the left hand side of (\ref{4.11}). We consider the boundary value problem
\begin{equation}\label{4.13}
\begin{array}{l}
  {\mathcal L}_2[\widetilde H]=\Phi_2(\varphi,\widetilde H)\quad\mbox{in}\ I, \\[0.2cm]
  \widetilde H(x_\ell)=0,\ \widetilde H(x_r)=0.
\end{array}
\end{equation}
The homogeneous boundary conditions appear naturally in the case of rotationally symmetric Willmore surfaces (see \cite{Mosel_01}, section 6.3).\\[0.1cm]
Note that ${\mathcal L}_2$ obeys the maximum principle if
\begin{equation}\label{4.14}
  \frac{f''}{f}-\frac{\alpha}{\beta}\,(1+f'^2)<0
  \quad\mbox{resp..}\quad
  \frac{1}{f^2}<\frac{\alpha}{\beta}\,.
\end{equation}
Let us work on this assumption.\\[0.1cm]
We remark that for sufficiently small perturbations $\|\psi\|_{2+\alpha,I}<\varepsilon$ die residual sums in (\ref{4.3}) and (\ref{4.12}) converge. For $\varepsilon<1$ and $\|\widetilde H\|_{2+\alpha,I}<1$ we find constants ${\mathcal C}_4,{\mathcal C}_5\in(0,+\infty)$ such that
\begin{equation}\label{4.15}
\begin{array}{lll}
  \|\Phi_1(\varphi,\widetilde H)\|_{\alpha,I}\negthickspace
  & \le & \negthickspace\displaystyle
          {\mathcal C}_4(\|\widetilde H\|_{2+\alpha,I}+\|\varphi\|_{2+\alpha,I}^2), \\[0.4cm]
  \|\Phi_2(\varphi,\widetilde H)\|_{\alpha,I}\negthickspace
  & \le & \negthickspace\displaystyle
          {\mathcal C}_4(\|\widetilde H\|_{2+\alpha,I}^3+\|\varphi\|_{2+\alpha,I}\|\widetilde H\|_{2+\alpha,I})
\end{array}
\end{equation}
as well as
\begin{equation}\label{4.16}
\begin{array}{l}
  \|\Phi_1(\varphi_1,\widetilde H_1)-\Phi_1(\varphi_2,\widetilde H_2)\|_{\alpha,I} \\[0.2cm]
  \hspace*{1.2cm}\displaystyle
  \le\,{\mathcal C}_5
       \{\|\widetilde H_1-\widetilde H_2\|_{2+\alpha,I}
         +(\|\varphi_1\|_{2+\alpha,I}+\|\varphi_2\|_{2+\alpha,I})\|\varphi_1-\varphi_2\|_{2+\alpha,I}\}, \\[0.2cm]
  \|\Phi_2(\varphi_1,\widetilde H_1)-\Phi_2(\varphi_2,\widetilde H_2)\|_{\alpha,I} \\[0.2cm]
  \hspace*{1.2cm}\displaystyle
  \le\,{\mathcal C}_5
       (\|\varphi_1\|_{2+\alpha,I}+\|\varphi_2\|_{2+\alpha,I}
        +\|\widetilde H_1\|_{2+\alpha,I}+\|\widetilde H_2\|_{2+\alpha,I})\cdot\ldots \\[0.2cm]
  \hspace*{2cm}\displaystyle
       \ldots\cdot(\|\varphi_1-\varphi_2\|_{2+\alpha,I}+\|\widetilde H_1-\widetilde H_2\|_{2+\alpha,I}).
\end{array}
\end{equation}
We will solve (\ref{4.9}) and (\ref{4.13}) successively: Let $f$ be stable, and let (\ref{4.14}) be fulfilled. We start with a pair $(\varphi_0,\widetilde H_0).$ First, with the right hand side $\Phi_2(\varphi_0,\widetilde H_0)$ we get an unique solution $\widetilde H_1$ from (\ref{4.13}). Now, we solve (\ref{4.9}) with $\Phi_1(\varphi_0,\widetilde H_1)$ to get an unique solution $\varphi_1.$ In this way we proceed with $(\varphi_1,\widetilde H_1).$\\[0.1cm]
Let $\|\varphi_0\|_{2+\alpha,I},\|\widetilde H_0\|_{2+\alpha,I}<\varepsilon$ with the above $\varepsilon\in(0,1).$ The Schauder estimates
\begin{equation}\label{4.17}
  \|\widetilde H_1\|_{2+\alpha,I}
  \le{\mathcal C}_6\|\Phi_2(\varphi_0,\widetilde H_0)\|_{\alpha,I}\,,\quad
  \|\varphi_1\|_{2+\alpha,I}
  \le{\mathcal C}_6\|\Phi_1(\varphi_0,\widetilde H_1)\|_{\alpha,I}
\end{equation}
with a constant ${\mathcal C}_6\in(0,+\infty)$ due to (\ref{6.40}), (\ref{4.7}) and (\ref{4.11}) yield (see (\ref{4.15}))
\begin{equation}\label{4.18}
\begin{array}{l}
  \|\widetilde H_1\|_{2+\alpha,I}
  \le{\mathcal C}_4{\mathcal C}_6(\varepsilon^3+\varepsilon^2)
  \le 2{\mathcal C}_4{\mathcal C}_6\varepsilon^2\,, \\[0.2cm]
  \|\varphi_1\|_{2+\alpha,I}
  \le{\mathcal C}_4{\mathcal C}_6(\|\widetilde H_1\|_{2+\alpha,I}+\varepsilon^2)
  \le{\mathcal C}_4{\mathcal C}_6(1+2{\mathcal C}_4{\mathcal C}_6)\varepsilon^2\,.
\end{array}
\end{equation}
We choose
\begin{equation}\label{4.19}
  \varepsilon
  \le\min
       \left\{
         1,
         \frac{1}{2{\mathcal C}_4{\mathcal C}_6}\,,
         \frac{1}{{\mathcal C}_4{\mathcal C_6}(1+{\mathcal C}_4{\mathcal C}_6)}
       \right\}.
\end{equation}
Then there follow $\|\widetilde H_1\|_{2+\alpha,I}\le\varepsilon$ and $\|\varphi_1\|_{2+\alpha,I}\le\varepsilon.$ We continue this procedure to get
\begin{equation}\label{4.20}
  \|\widetilde H_k\|_{2+\alpha,I}\le\varepsilon,\quad
  \|\varphi_k\|_{2+\alpha,I}\le\varepsilon
  \quad\mbox{for}\ k=0,1,2,\ldots
\end{equation}
From (\ref{4.16}) we deduce
\begin{equation}\label{4.21}
  \|\widetilde H_{k+1}-\widetilde H_k\|_{2+\alpha,I}
  \le 4\varepsilon\,{\mathcal C}_5{\mathcal C}_6
      (\|\widetilde H_k-\widetilde H_{k-1}\|_{2+\alpha,I}
       +\|\varphi_k-\varphi_{k-1}\|_{2+\alpha,I})
\end{equation}
as well as (for ${\mathcal C}_5$ and ${\mathcal C}_6$ be sufficiently large)
\begin{equation}\label{4.22}
\begin{array}{lll}
  \|\varphi_{k+1}-\varphi_k\|_{2+\alpha,I}\negthickspace
  & \le & \negthickspace\displaystyle
          {\mathcal C}_5{\mathcal C}_6\|\widetilde H_{k+1}-\widetilde H_k\|_{2+\alpha,I}
          +2\varepsilon\,{\mathcal C}_5{\mathcal C}_6\|\varphi_k-\varphi_{k-1}\|_{2+\alpha,I} \\[0.2cm]
  & \le & \negthickspace\displaystyle
          4\varepsilon\,{\mathcal C}_5^2{\mathcal C}_6^2\|\widetilde H_k-\widetilde H_{k-1}\|_{2+\alpha,I}
          +6\varepsilon\,{\mathcal C}_5^2{\mathcal C}_6^2\|\varphi_k-\varphi_{k-1}\|_{2+\alpha,I}\,.
\end{array}
\end{equation}
From the last inequalities we conclude
\begin{equation}\label{4.23}
\begin{array}{l}
  \|\varphi_{k+1}-\varphi_k\|_{2+\alpha,I}
  +\|\widetilde H_{k+1}-\widetilde H_k\|_{2+\alpha,I} \\[0.2cm]
  \hspace*{1.2cm}\displaystyle
  \le 10\varepsilon\,{\mathcal C}_5^2{\mathcal C}_6^2
      (\|\varphi_k-\varphi_{k-1}\|_{2+\alpha,I}
       +\|\widetilde H_k-\widetilde H_{k-1}\|_{2+\alpha,I})
\end{array}
\end{equation}
for $k=1,2,\ldots$ Additionally we assume the smallness property
\begin{equation}\label{4.24}
  \varepsilon<\frac{1}{20\,{\mathcal C}_5^2{\mathcal C}_6^2}
\end{equation}
such that the sequence $(\varphi_k,\widetilde H_k)_{k=1,2,\ldots}$ is contractive. Then we find solutions $\varphi\in C^{2+\alpha}(I,\mathbb R)$ and $\widetilde H\in C^{2+\alpha}(I,\mathbb R),$ $\|\varphi\|_{2+\alpha,I}<\varepsilon$ and $\|\widetilde H\|_{2+\alpha,I}<\varepsilon,$ of the non-linear coupled boundary value problems (\ref{4.9}) and (\ref{4.13}). As a meridian curve of the rotationally symmetric immersion with mean curvature $H\in C^{2+\alpha}(I,\mathbb R),$ we have $f\in C^{4+\alpha}(I,\mathbb R)$ due to well-known regularity results.
\begin{theorem}
Let $f\in C^{4+\alpha}(I,\mathbb R)$ be a stable solution of the minimal surface equation (\ref{2.2}) for given boundary values $f_\ell$ and $f_r.$ Let $f$ satisfy (\ref{4.14}). Then we find a solution pair $(f,H)$ of (\ref{3.16}) and (\ref{3.19}) as a sufficiently small perturbation in the $C^{4+\alpha}$-norm of the minimal surface.
\end{theorem}
\subsection{Numerical results}
\setcounter{equation}{0}
\setcounter{Lem}{0}
\setcounter{Prop}{0}
\setcounter{Def}{0}
\setcounter{Theo}{0}
\setcounter{Cor}{0}
\subsubsection{The catenoid}
Two coaxial circular boundary curves of common radius $r>0$ span a minimal surface if the distance $h>0$ between the rings is sufficiently small.\\[0.1cm]
We investigate the dependence between the ratio $\nicefrac{h}{r}$ and the surface area $A$ numerically. For this, we increase the distance $h$ of the rings starting from $h=0.1$ by adding successively $\triangle h=0.1$ while keeping the radius $r=1.5088795$ fixed (all calculations were done with Ken Brakke's Surface Evolver \cite{Brakke_01}).\\[0.1cm]
In the case $\nicefrac{h}{r}=1.3256$ there is no two-fold connected minimal surface. The given surface area is the total area of the two discs of radius $r.$\\[0.6cm]
\small{
\begin{tabular}{c|cccccccccc}
  $\bef{\nicefrac{h}{r}}$ & $0.0663$ & $0.1325$ & $0.1988$ & $0.2651$ & $0.3314$
                          & $0.3976$ & $0.4639$ & $0.5302$ & $0.5965$ & $0.6627$ \\
  \hline
  $\bef{A}$               & $0.9053$ & $1.8881$ & $2.8292$ & $3.7667$ & $4.7001$ 
                          & $5.6285$ & $6.5506$ & $7.4656$ & $8.3723$ & $9.2694$
\end{tabular}\\[0.4cm]
\begin{tabular}{c|cccccccccc}
  $\bef{\nicefrac{h}{r}}$ & $0.7290$  & $0.7953$  & $0.8616$  & $0.9278$  & $0.9941$
                          & $1.0604$  & $1.1267$  & $1.1929$  & $1.2592$  & $1.3256$ \\
  \hline
  $\bef{A}$               & $10.1558$ & $11.0302$ & $11.8810$ & $12.7189$ & $13.5347$
                          & $14.3377$ & $15.0976$ & $15.8227$ & $16.5026$ & $14.3250$
\end{tabular}
}
\vspace*{0.2cm}
\begin{center}
\pspicture(0,0)(12.6,2.6)
\rput(1,1.5){\fbox{\epsfig{file=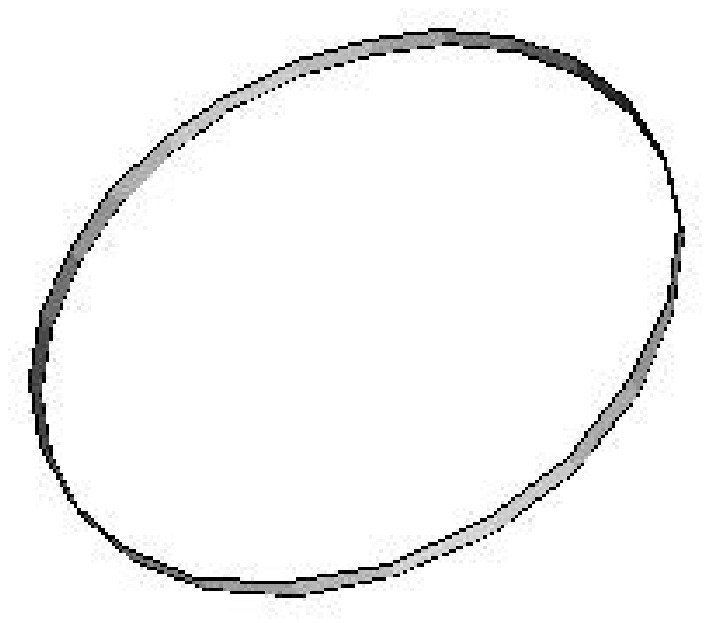, height=2.2cm}}}
\rput(1,0){\scalebox{0.8}{$\nicefrac{h}{r}=0.0663$}}
\rput(4.5,1.5){\fbox{\epsfig{file=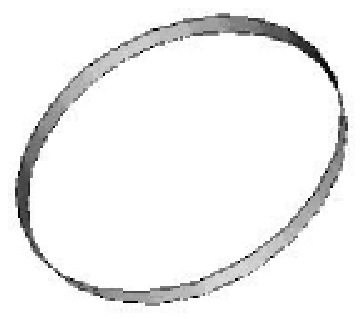, height=2.2cm}}}
\rput(4.5,0){\scalebox{0.8}{$\nicefrac{h}{r}=0.1325$}}
\rput(8,1.5){\fbox{\epsfig{file=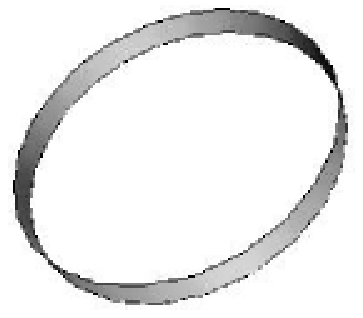, height=2.2cm}}}
\rput(8,0){\scalebox{0.8}{$\nicefrac{h}{r}=0.1988$}}
\rput(11.5,1.5){\fbox{\epsfig{file=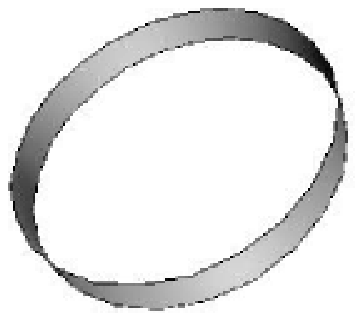, height=2.2cm}}}
\rput(11.5,0){\scalebox{0.8}{$\nicefrac{h}{r}=0.2651$}}
\endpspicture
\end{center}
\vspace*{0.1cm}
\begin{center}
\pspicture(0,0)(12.6,2.6)
\rput(1,1.5){\fbox{\epsfig{file=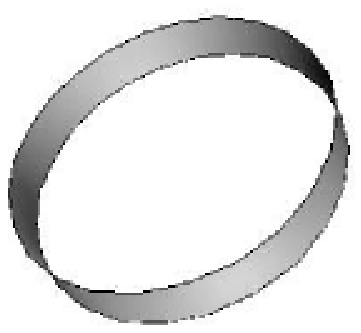, height=2.2cm}}}
\rput(1,0){\scalebox{0.8}{$\nicefrac{h}{r}=0.3314$}}
\rput(4.5,1.5){\fbox{\epsfig{file=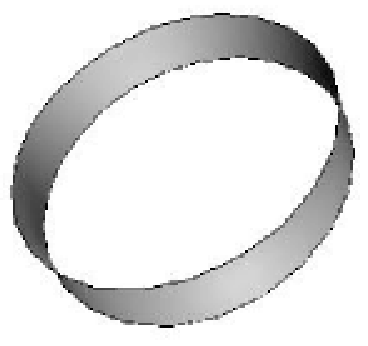, height=2.2cm}}}
\rput(4.5,0){\scalebox{0.8}{$\nicefrac{h}{r}=0.3976$}}
\rput(8,1.5){\fbox{\epsfig{file=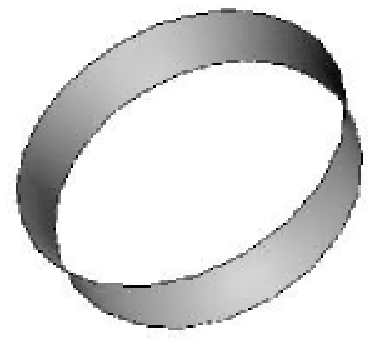, height=2.2cm}}}
\rput(8,0){\scalebox{0.8}{$\nicefrac{h}{r}=0.4639$}}
\rput(11.5,1.5){\fbox{\epsfig{file=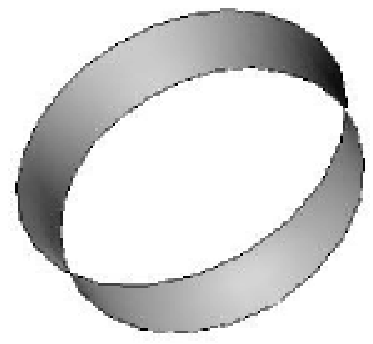, height=2.2cm}}}
\rput(11.5,0){\scalebox{0.8}{$\nicefrac{h}{r}=0.5302$}}
\endpspicture
\end{center}
\vspace*{0.1cm}
\begin{center}
\pspicture(0,0)(12.6,2.6)
\rput(1,1.5){\fbox{\epsfig{file=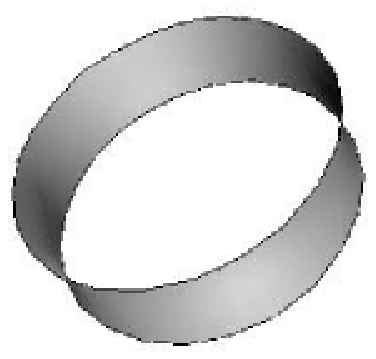, height=2.2cm}}}
\rput(1,0){\scalebox{0.8}{$\nicefrac{h}{r}=0.5965$}}
\rput(4.5,1.5){\fbox{\epsfig{file=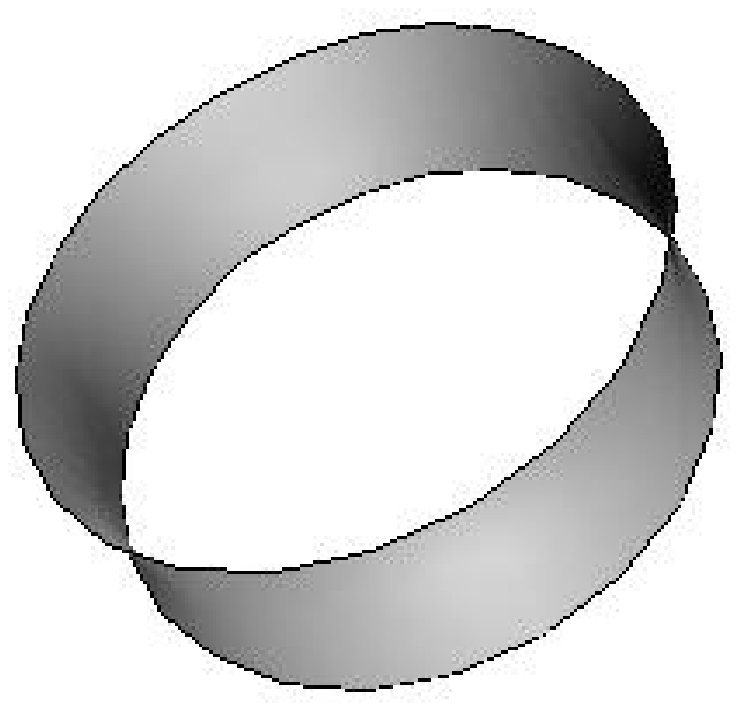, height=2.2cm}}}
\rput(4.5,0){\scalebox{0.8}{$\nicefrac{h}{r}=0.6627$}}
\rput(8,1.5){\fbox{\epsfig{file=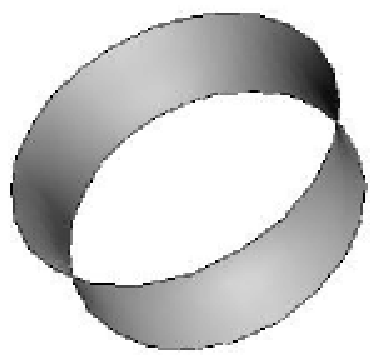, height=2.2cm}}}
\rput(8,0){\scalebox{0.8}{$\nicefrac{h}{r}=0.7290$}}
\rput(11.5,1.5){\fbox{\epsfig{file=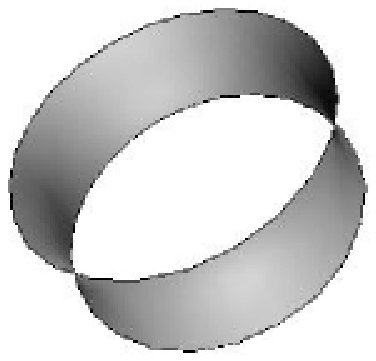, height=2.2cm}}}
\rput(11.5,0){\scalebox{0.8}{$\nicefrac{h}{r}=0.7953$}}
\endpspicture
\end{center}
\vspace*{0.1cm}
\begin{center}
\pspicture(0,0)(12.6,2.6)
\rput(1,1.5){\fbox{\epsfig{file=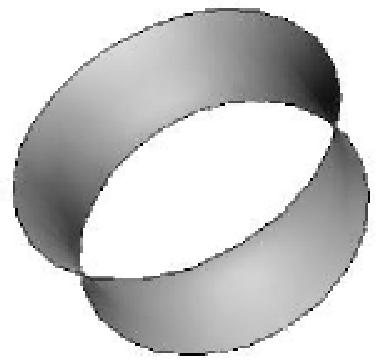, height=2.2cm}}}
\rput(1,0){\scalebox{0.8}{$\nicefrac{h}{r}=0.8616$}}
\rput(4.5,1.5){\fbox{\epsfig{file=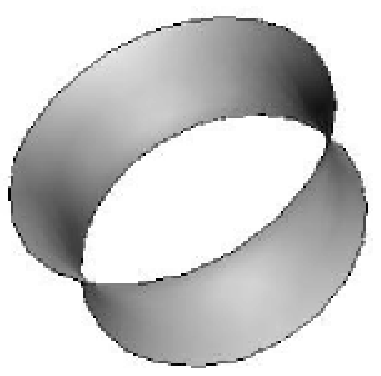, height=2.2cm}}}
\rput(4.5,0){\scalebox{0.8}{$\nicefrac{h}{r}=0.9278$}}
\rput(8,1.5){\fbox{\epsfig{file=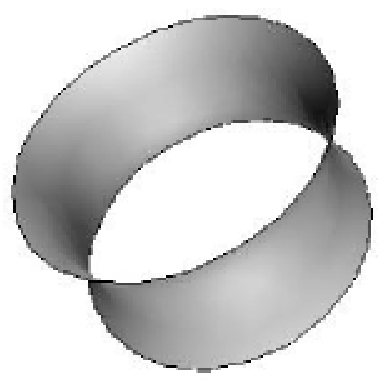, height=2.2cm}}}
\rput(8,0){\scalebox{0.8}{$\nicefrac{h}{r}=0.9941$}}
\rput(11.5,1.5){\fbox{\epsfig{file=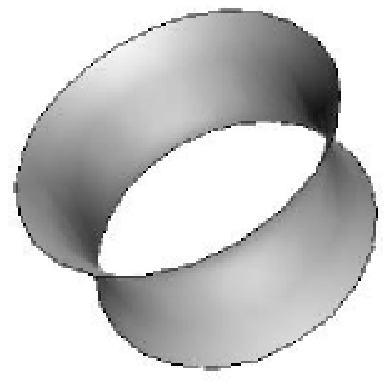, height=2.2cm}}}
\rput(11.5,0){\scalebox{0.8}{$\nicefrac{h}{r}=1.0604$}}
\endpspicture
\end{center}
\vspace*{0.1cm}
\begin{center}
\pspicture(0,0)(12.6,2.6)
\rput(1,1.5){\fbox{\epsfig{file=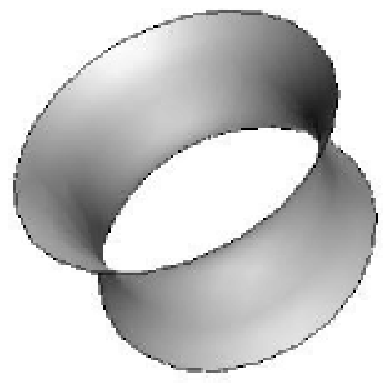, height=2.2cm}}}
\rput(1,0){\scalebox{0.8}{$\nicefrac{h}{r}=1.1267$}}
\rput(4.5,1.5){\fbox{\epsfig{file=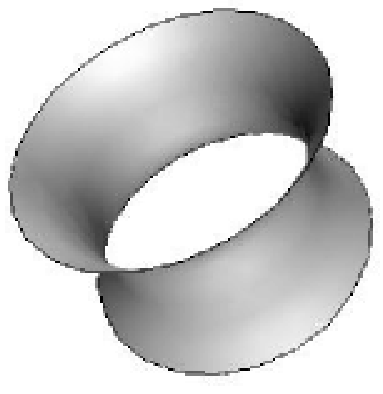, height=2.2cm}}}
\rput(4.5,0){\scalebox{0.8}{$\nicefrac{h}{r}=1.1929$}}
\rput(8,1.5){\fbox{\epsfig{file=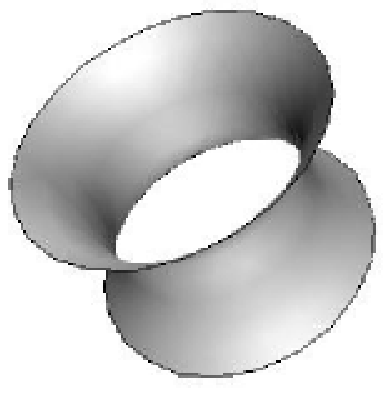, height=2.2cm}}}
\rput(8,0){\scalebox{0.8}{$\nicefrac{h}{r}=1.2592$}}
\rput(11.5,1.5){\fbox{\epsfig{file=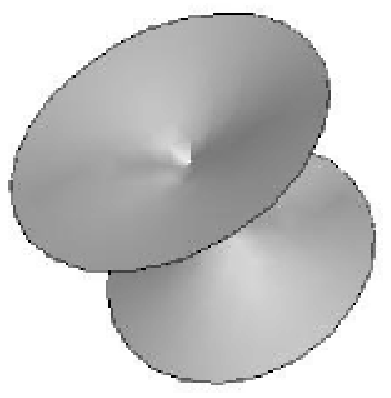, height=2.2cm}}}
\rput(11.5,0){\scalebox{0.8}{$\nicefrac{h}{r}=1.3256$}}
\endpspicture
\end{center}
\subsubsection{The Willmore catenoid}
We investigate numerically critical points of Willmore's functional (\ref{1.4}) with $\alpha=0,$ $\gamma=0,$ $\beta=1.$ The following tables compare the surface area $A$ with the Willmore energy $E$ for different distances $h$ ($r=1$).\\[0.1cm]
\small{
\begin{tabular}{c|ccccccccccc}
  $\bef{h}$ & $1.0$  & $1.1$  & $1.2$  & $1.3$  & $1.4$  & $1.5$  & $1.6$  & $1.7$  & $1.8$  & $1.9$  & $2.0$   \\
  \hline
  $\bef{A}$ & $5.98$ & $6.50$ & $6.98$ & $7.42$ & $7.77$ & $8.13$ & $8.49$ & $8.87$ & $9.27$ & $9.68$ & $10.11$ \\
  \hline
  $\bef{E}$ & $0.00$ & $0.00$ & $0.00$ & $0.00$ & $0.02$ & $0.10$ & $0.22$ & $0.38$ & $0.56$ & $0.77$ & $0.98$
\end{tabular}\\[0.3cm]
\begin{tabular}{c|ccccccccc}
  $\bef{h}$ & $2.1$   & $2.2$   & $2.3$   & $2.4$   & $2.5$   & $2.6$   & $2.7$   & $2.8$   & $2.9$   \\
  \hline
  $\bef{A}$ & $10.56$ & $11.02$ & $11.51$ & $12.01$ & $12.54$ & $13.10$ & $13.66$ & $14.25$ & $14.86$ \\
  \hline
  $\bef{E}$ & $1.20$  & $1.43$  & $1.67$  & $1.90$  & $2.13$  & $2.36$  & $2.58$  & $2.80$  & $3.02$
\end{tabular}\\[0.3cm]
\begin{tabular}{c|cccccccc}
  $\bef{h}$ & $3.0$   & $4.0$   & $5.0$   & $6.0$   & $7.0$   & $8.0$   & $9.0$    & $10.0$   \\
  \hline
  $\bef{A}$ & $15.50$ & $23.21$ & $34.15$ & $48.08$ & $65.37$ & $85.79$ & $111.07$ & $145.22$ \\
  \hline
  $\bef{E}$ & $3.23$  & $5.02$  & $6.62$  & $7.56$  & $8.25$  & $8.72$  & $9.26$   & $9.71$
\end{tabular}
}
\begin{center}
\pspicture(0,0)(12.6,2.6)
\rput(1,1.5){\fbox{\epsfig{file=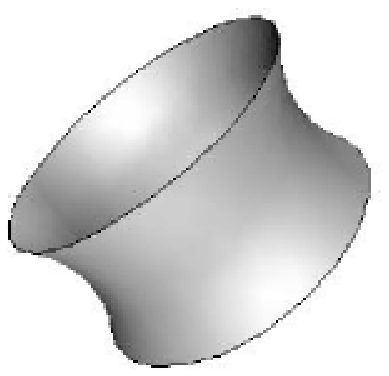, height=2.2cm}}}
\rput(1,0){\scalebox{0.8}{$h=10$}}
\rput(4.5,1.5){\fbox{\epsfig{file=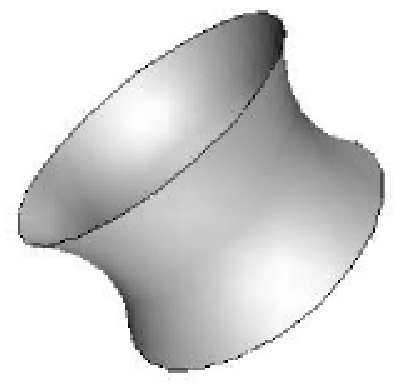, height=2.2cm}}}
\rput(4.5,0){\scalebox{0.8}{$h=11$}}
\rput(8,1.5){\fbox{\epsfig{file=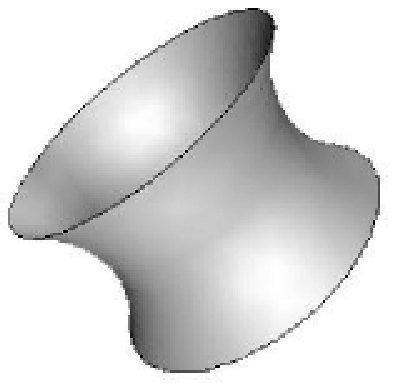, height=2.2cm}}}
\rput(8,0){\scalebox{0.8}{$h=12$}}
\rput(11.5,1.5){\fbox{\epsfig{file=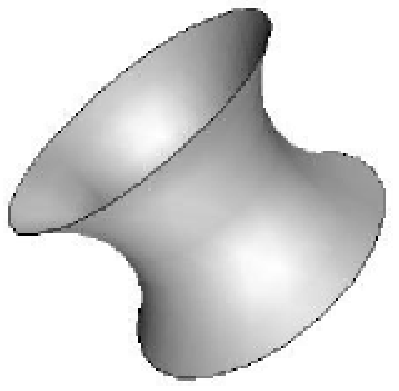, height=2.2cm}}}
\rput(11.5,0){\scalebox{0.8}{$h=13$}}
\endpspicture
\end{center}
\vspace*{0.1cm}
\begin{center}
\pspicture(0,0)(12.6,2.6)
\rput(1,1.5){\fbox{\epsfig{file=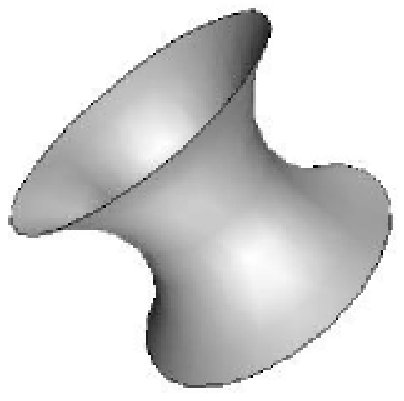, height=2.2cm}}}
\rput(1,0){\scalebox{0.8}{$h=14$}}
\rput(4.5,1.5){\fbox{\epsfig{file=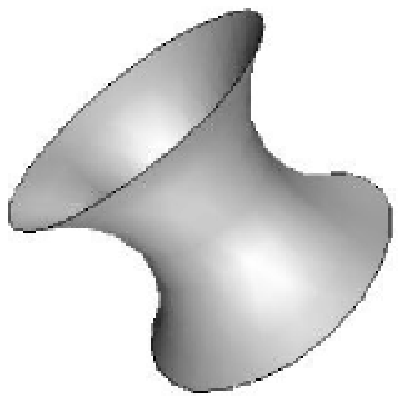, height=2.2cm}}}
\rput(4.5,0){\scalebox{0.8}{$h=15$}}
\rput(8,1.5){\fbox{\epsfig{file=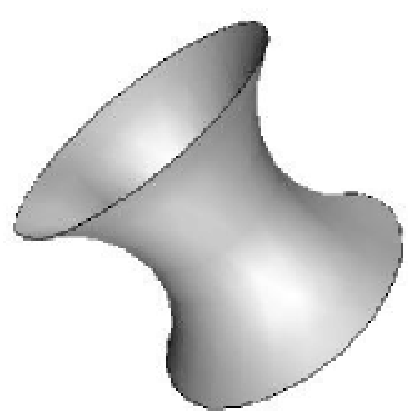, height=2.2cm}}}
\rput(8,0){\scalebox{0.8}{$h=16$}}
\rput(11.5,1.5){\fbox{\epsfig{file=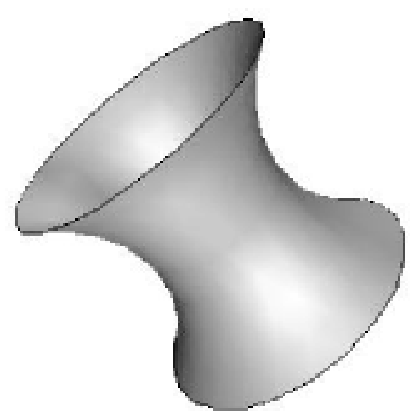, height=2.2cm}}}
\rput(11.5,0){\scalebox{0.8}{$h=17$}}
\endpspicture
\end{center}
\vspace*{0.1cm}
\begin{center}
\pspicture(0,0)(12.6,2.6)
\rput(1,1.5){\fbox{\epsfig{file=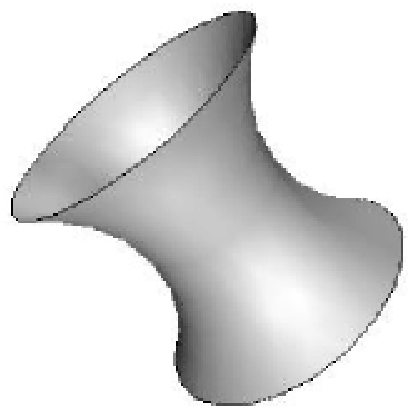, height=2.2cm}}}
\rput(1,0){\scalebox{0.8}{$h=18$}}
\rput(4.5,1.5){\fbox{\epsfig{file=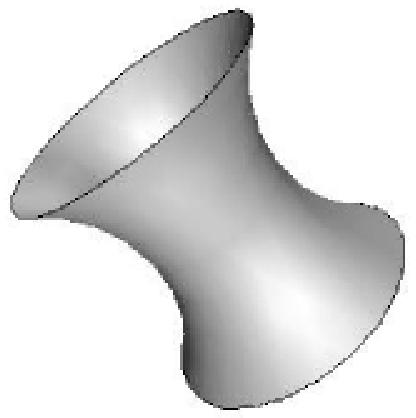, height=2.2cm}}}
\rput(4.5,0){\scalebox{0.8}{$h=19$}}
\rput(8,1.5){\fbox{\epsfig{file=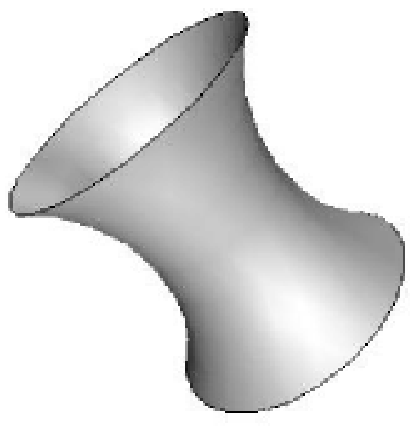, height=2.2cm}}}
\rput(8,0){\scalebox{0.8}{$h=20$}}
\rput(11.5,1.5){\fbox{\epsfig{file=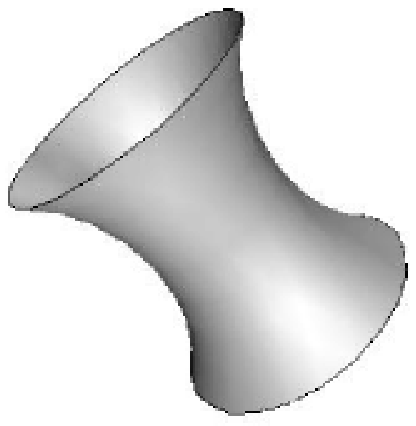, height=2.2cm}}}
\rput(11.5,0){\scalebox{0.8}{$h=21$}}
\endpspicture
\end{center}
\vspace*{0.1cm}
\begin{center}
\pspicture(0,0)(12.6,2.6)
\rput(1,1.5){\fbox{\epsfig{file=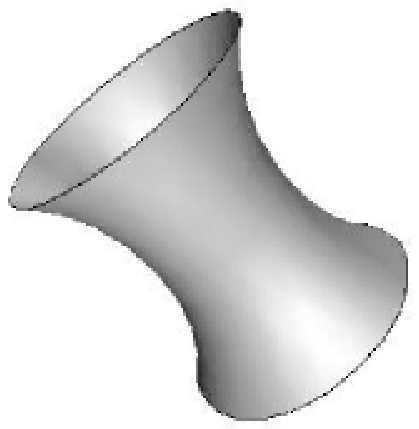, height=2.2cm}}}
\rput(1,0){\scalebox{0.8}{$h=22$}}
\rput(4.5,1.5){\fbox{\epsfig{file=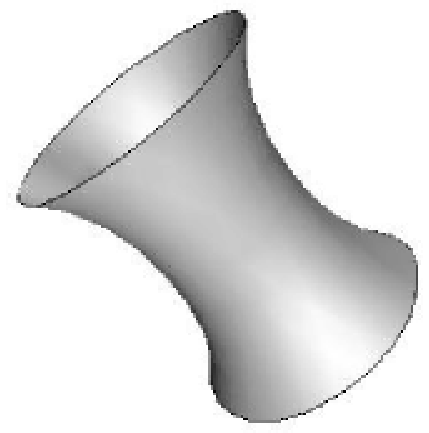, height=2.2cm}}}
\rput(4.5,0){\scalebox{0.8}{$h=23$}}
\rput(8,1.5){\fbox{\epsfig{file=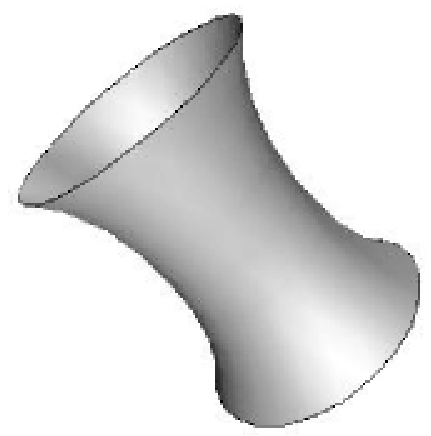, height=2.2cm}}}
\rput(8,0){\scalebox{0.8}{$h=24$}}
\rput(11.5,1.5){\fbox{\epsfig{file=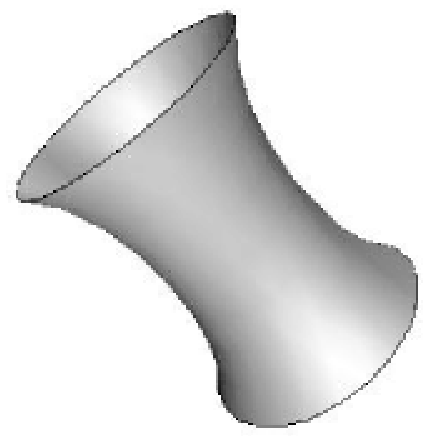, height=2.2cm}}}
\rput(11.5,0){\scalebox{0.8}{$h=25$}}
\endpspicture
\end{center}
\vspace*{0.1cm}
\begin{center}
\pspicture(0,0)(12.6,2.6)
\rput(1,1.5){\fbox{\epsfig{file=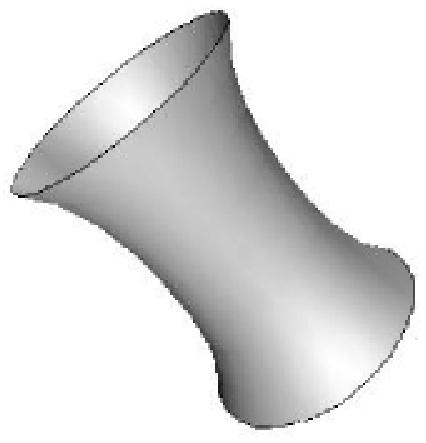, height=2.2cm}}}
\rput(1,0){\scalebox{0.8}{$h=26$}}
\rput(4.5,1.5){\fbox{\epsfig{file=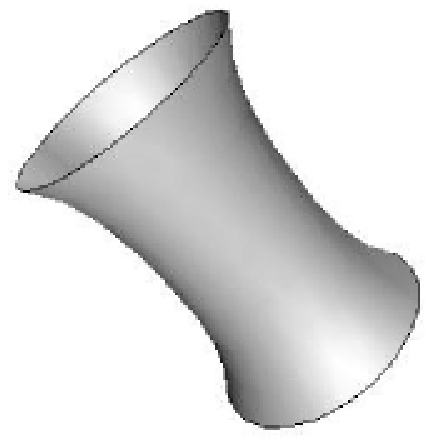, height=2.2cm}}}
\rput(4.5,0){\scalebox{0.8}{$h=27$}}
\rput(8,1.5){\fbox{\epsfig{file=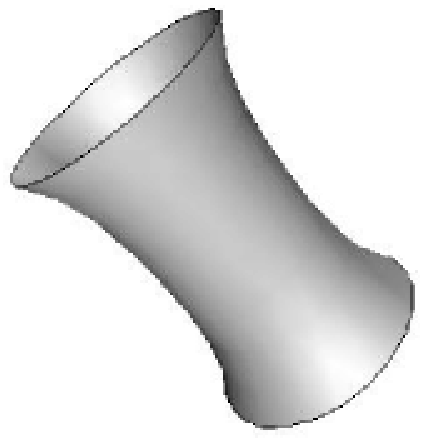, height=2.2cm}}}
\rput(8,0){\scalebox{0.8}{$h=28$}}
\rput(11.5,1.5){\fbox{\epsfig{file=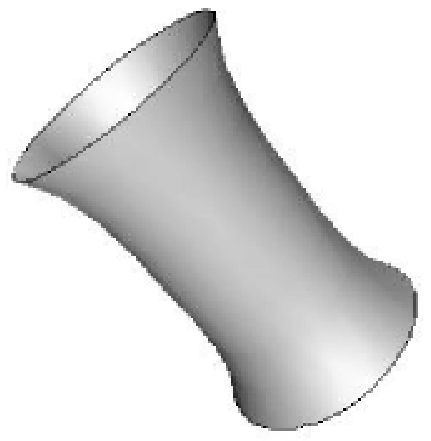, height=2.2cm}}}
\rput(11.5,0){\scalebox{0.8}{$h=29$}}
\endpspicture
\end{center}
It is conjectured that such catenoid-type critical points of Willmore's functional exist for all distances $h>0.$ Furthermore, for $h$ sufficiently large the critical immersions seem to tend to the sphere $S^2$ - except from a small neighborhood of the boundary curves where the mean curvature $H$ vanishes. Therefore, we expect ${\mathcal W}\to 4\pi\approx 12.56$ for $r\to\infty$ due to the fact that $4\pi$ is exactly the Willmore energy for the round sphere.
\begin{center}
\pspicture(0,0)(12.6,2.6)
\rput(1,1.5){\fbox{\epsfig{file=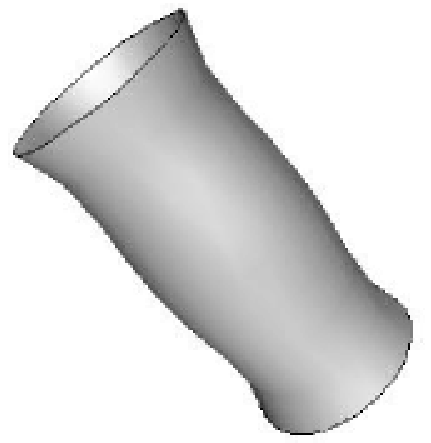, height=2.2cm}}}
\rput(1,0){\scalebox{0.8}{$h=40$}}
\rput(4.5,1.5){\fbox{\epsfig{file=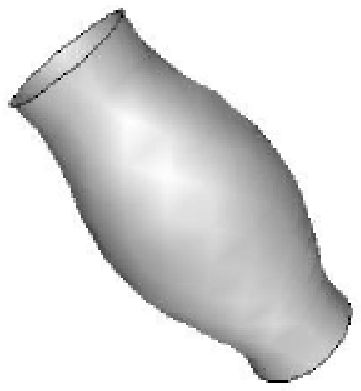, height=2.2cm}}}
\rput(4.5,0){\scalebox{0.8}{$h=60$}}
\rput(8,1.5){\fbox{\epsfig{file=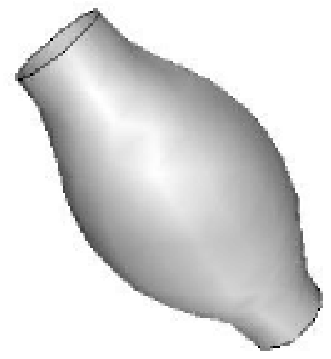, height=2.2cm}}}
\rput(8,0){\scalebox{0.8}{$h=80$}}
\rput(11.5,1.5){\fbox{\epsfig{file=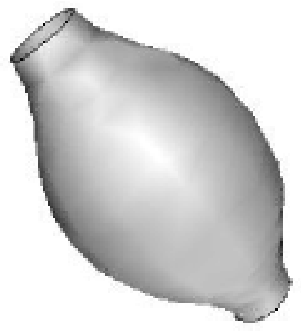, height=2.2cm}}}
\rput(11.5,0){\scalebox{0.8}{$h=100$}}
\endpspicture
\end{center}
\subsection{Appendix: Linear ordinary differential equations}
\subsubsection{The maximum principle}
\setcounter{equation}{0}
\setcounter{Lem}{0}
\setcounter{Prop}{0}
\setcounter{Def}{0}
\setcounter{Theo}{0}
\setcounter{Cor}{0}
Let $u\in C^{2+\alpha}(I,\mathbb R),$ $\alpha\in(0,1),$ $I=[a,b],$ solve the boundary value problem
  $$\widetilde{\mathcal L}[u](x)
    :=u''(x)+\widetilde p(x)u'(x)+\widetilde q(x)\widetilde u(x)=\widetilde f(x)\quad\mbox{in}\ I,
    \quad u(a)=\eta_1\,,\ u(b)=\eta_2\,,\eqno(\mbox{6.1a})$$
where $\widetilde p,\widetilde q\in C^\alpha(I,\mathbb R)$ and $f\in C^\alpha(I,\mathbb R).$ Multiplication with the positive function
\setcounter{equation}{1}
\begin{equation}\label{6.2}
  p(x):=\exp\int\widetilde p(x)\,dx
\end{equation}
gives the Sturm-Liouville boundary value problem
  $${\mathcal L}[u]:=(pu')'+qu=f\quad\mbox{in}\ I,
    \quad u(a)=\eta_1\,,\ u(b)=\eta_2\,,\eqno(\mbox{6.1b})$$
where $q:=p\widetilde q$ and $f:=p\widetilde f.$ Note that $p\in C^1(I,\mathbb R).$\\[0.1cm]
In \cite{Walter_01}, \S 26, Satz III we find the fundamental existence theorem:
\begin{proposition}
(Existence and uniqueness)\\[0.1cm]
The inhomogeneous problem (6.1b) has a unique solution $u\in C^{2+\alpha}(I,\mathbb R)$ if and only if the homogeneous problem
\begin{equation}\label{6.3}
  {\mathcal L}[u]=0,\quad u(a)=0,\ u(b)=0,
\end{equation}
has only the trivial solution $u\equiv 0.$
\end{proposition}
\noindent
The H\"older regularity for $u$ comes from the assumption $f\in C^\alpha(I,\mathbb R).$ Furthermore, the uniqueness has its origin in the following maximum principle which can be found in \cite{Walter_01}, Erg\"anzung I for \S 26.
\begin{proposition}
(Maximum principle)\\[0.1cm]
Let $u\in C^2(I,\mathbb R)$ solve (6.1b) with $q\le 0$ (which is equivalent to $\widetilde q\le 0$).
\begin{itemize}
\item[(i)]
If ${\mathcal L}[u]\ge 0$ and $u$ has an inner positive maximum, then $u\equiv\mbox{const.}$ 
\item[(ii)]
If ${\mathcal L}[u]\le 0$ and $u$ has an inner negative minimum, then $u\equiv\mbox{const.}$
\end{itemize}
\end{proposition}
\begin{corollary}
(Uniqueness)\\[0.1cm]
Let $u,v\in C^2(I,\mathbb R)$ solve
\begin{equation}\label{6.4}
\begin{array}{l}
  {\mathcal L}[u]=f,\quad u(a)=\eta_1\,,\ u(b)=\eta_2\,, \\[0.2cm]
  {\mathcal L}[v]=f,\quad v(a)=\eta_1\,,\ v(b)=\eta_2\,.
\end{array}
\end{equation}
Assume that $q\le 0.$ Then it holds $u\equiv v$ in $I.$
\end{corollary}
\begin{proof}
Consider $w:=u-v$ with ${\mathcal L}[w]=0$ and $w(a)=0,$ $w(b)=0.$ The maximum principle yields $w\equiv 0.$
\end{proof}
\noindent
For $q\le 0$ the homogeneous problem has only the trivial solution, that is (6.1b) is uniquely solvable.\\[0.1cm]
Alternatively, if the homogeneous problem has non-trivial solutions, the inhomogeneous problem may have no solutions as the following example shows (\cite{Walter_01}, \S 26):
\begin{equation}\label{6.5}
  u''(x)+u(x)=1,\quad u(0)=0,\ u(\pi)=1.
\end{equation}
\subsubsection{A priori estimates}
Let $p_0,p_1\in\mathbb R$ and $p_1'\in\mathbb R$ be constants such that
\begin{equation}\label{6.6}
  0<p_0\le p(x)\le p_1<+\infty
  \quad\mbox{and}\quad
  |p'(x)|\le p_1'<+\infty
  \quad\mbox{for all}\ x\in I.
\end{equation}
The maximum principle from Proposition 6.2 gives the following bounds on the $C^0$-norm.
\begin{corollary}
(Estimate of the $C^0$-norm)\\[0.1cm]
Let $u\in C^2(I,\mathbb R)$ solve (6.1b) with $q\le 0.$ Then
\begin{equation}\label{6.7}
  \|u\|_{0,I}
  \le 3\max_{x\in\partial I}|u(x)|
      +\frac{\|f\|_{0,I}}{\nu e^{\nu a}}\,e^{\nu(b-a)}\,,
  \quad\nu:=\frac{1+p_1'}{p_0}\,.
\end{equation}
\end{corollary}
\begin{proof}
\begin{itemize}
\item[1.]
Consider the positive functions
\begin{equation}\label{6.8}
  v(x):=\max_{x\in\partial I}|u(x)|+\frac{\|f\|_{0,I}}{\mu e^{\mu a}}\,(e^{b\mu}-e^{\mu x})\,,\quad x\in I,
\end{equation}
where $\mu p_0-p_1'=1$ (note that $\mu>0$), such that $p'+p\mu\ge 1.$ It follows that
\begin{equation}\label{6.9}
  [p(e^{\mu x})']'=\mu[pe^{\mu x}]'=\mu[p'+p\mu]e^{\mu x}\ge\mu e^{\mu x}\ge\mu e^{\mu a}\,.
\end{equation}
Because $qv\le 0$ we conclude
\begin{equation}\label{6.10}
  {\mathcal L}[v]
  =(pv')'+qv
  \le -\frac{\|f\|_{0,I}}{\mu e^{\mu a}}\,\mu e^{\mu a}+qv
  \le -\|f\|_{0,I}\,.
\end{equation}
Thus, ${\mathcal L}[v-u]\le 0.$ Due to the maximum principle, $v-u$ has no local negative minimum if it is not constant. Because $(v-u)|_{\partial I}\ge 0$ we have $v\ge u$ if $v-u\not\equiv\mbox{const}.$
\item[2.]
Let $\widetilde v=-v,$ then ${\mathcal L}[\widetilde v]\ge\|f\|_{0,I}.$ We have ${\mathcal L}[\widetilde v-u]\ge 0.$ Due to the maximum principle, $\widetilde v-u$ has no local positive maximum if it is not constant. Because $(\widetilde v-u)|_{\partial I}\le 0$ we have $\widetilde v\le u$ if $\widetilde v-u\not\equiv\mbox{const}.$
\item[3.]
If $v-u\not\equiv\mbox{const}$ and $\widetilde v-u\not\equiv\mbox{const}$ we have
\begin{equation}\label{6.11}
  -\max_{x\in\partial I}|u(x)|-\frac{\|f\|_{0,I}}{\mu e^{\mu a}}\,(e^{b\mu}-e^{\mu x})
  \le u(x)
  \le\max_{x\in\partial I}|u(x)|+\frac{\|f\|_{0,I}}{\mu e^{\mu a}}\,(e^{b\mu}-e^{\mu x})
  \quad\mbox{in}\ I.
\end{equation}
The statement follows for this case.
\item[4.]
Let $v-u\equiv C,$ $C\in\mathbb R.$ With $C=v(b)-u(b)=v(b)-\eta_2$ we get
\begin{equation}\label{6.12}
  |C|\le v(b)+|\eta_2|\le 2\max_{x\in\partial I}|u(x)|,
\end{equation}
and the statement follows from $|u|\le|v|+|C|.$ Analogously, for $\widetilde v-u\equiv{\widetilde C}$ we get $\displaystyle|\widetilde C|\le 2\max_{x\in\partial I}|u(x)|.$ This makes the proof complete.
\end{itemize}
\end{proof}
\noindent
In the following calculations we can replace $\|u\|_{0,I}$ by the above estimate if $q\le 0.$
\begin{lemma}
(Estimate of the $C^1$-norm)\\[0.1cm]
Let $u\in C^{2+\alpha}(I,\mathbb R)$ solve (6.1b). Then, for all real $0<\mu<\frac{1}{2}$ it holds
\begin{equation}\label{6.13}
  \|u\|_{1,I}
  \le\frac{b-a}{2p_0}\,\|f\|_{0,I}
     +\left(1+\frac{2p_1}{\mu p_0(b-a)}+\frac{q_1(b-a)}{2p_0}\right)\|u\|_{0,I}
     +\frac{\mu p_1(b-a)}{2p_0}\,\|u\|_{2+\alpha,I}\,.
\end{equation}
\end{lemma}
\begin{proof}
\begin{itemize}
\item[1.]
First, we estimate $|u'(c)|$ with $c:=\frac{a+b}{2}.$ Let $d:=\mu\frac{b-a}{2}$ with $\mu<\frac{1}{2}$ sufficiently small. Set $x_1:=c-d$ and $x_2:=c+d.$ By the mean value theorem there exists $\widetilde x\in[x_1,x_2]$ such that
\begin{equation}\label{6.14}
  |u'(\widetilde x)|
  =\frac{|u(x_1)-u(x_2)|}{|x_1-x_2|}
  =\frac{|u(x_1)-u(x_2)|}{2d}
  \le\frac{1}{d}\,\|u\|_{0,I}\,.
\end{equation}
Because $|c-\widetilde x|\le d$ we get
\begin{equation}\label{6.15}
\begin{array}{lll}
  |u'(c)|\negthickspace\displaystyle
  &  =  & \negthickspace\displaystyle
          |u'(\widetilde x)+\int\limits_{\widetilde x}^cu''(\xi)\,d\xi\,|
          \,\le\,\frac{1}{d}\,\|u\|_{0,I}+d\sup_{y\in[x_1,x_2]}|u''(y)| \\[0.8cm]
  & \le & \negthickspace\displaystyle
          \frac{2}{\mu(b-a)}\,\|u\|_{0,I}+\frac{\mu(b-a)}{2}\,[u]_{2,I}\,.
\end{array}
\end{equation}
\item[2.]
Integrating $(pu')'=f-qu$ from $\xi=c$ to $\xi=x$ yields
\begin{equation}\label{6.16}
  p(x)u'(x)=\int\limits_c^xf(\xi)\,d\xi-\int\limits_c^xq(\xi)u(\xi)\,d\xi+p(c)u'(c),
\end{equation}
that is,
\begin{equation}\label{6.17}
  [u]_{1,I}
  \le\frac{b-a}{2p_0}\,(\|f\|_{0,I}+q_1\|u\|_{0,I})
      +\frac{2p_1}{\mu p_0(b-a)}\,\|u\|_{0,I}
      +\frac{\mu p_1(b-a)}{2p_0}\,\|u\|_{2+\alpha,I}\,.
\end{equation}
The statement follows.
\end{itemize}
\end{proof}
\begin{lemma}
(Estimate of the $C^\alpha$-norm)\\[0.1cm]
Let $u\in C^{2+\alpha}(I,\mathbb R)$ solve (6.1b). Then
\begin{equation}\label{6.18}
  \|u\|_{\alpha,I}\le\|u\|_{0,I}+(b-a)^{1-\alpha}\|u\|_{1,I}
\end{equation}
with the above bounds on $\|u\|_{0,I}$ and $\|u\|_{1,I}.$
\end{lemma}
\begin{proof}
For arbitrary $x_1,x_2\in[a,b],$ $x_1\not=x_2,$ we have
\begin{equation}\label{6.19}
  |u(x_1)-u(x_2)|=|u'(\xi)||x_1-x_2|
\end{equation}
with $\xi\in[x_1,x_2].$ We conclude
\begin{equation}\label{6.20}
  \frac{|u(x_1)-u(x_2)|}{|x_1-x_2|^\alpha}
  =|u'(\xi)||x_1-x_2|^{1-\alpha}
  \le(b-a)^{1-\alpha}\|u\|_{1,I}\,,
\end{equation}
therefore
\begin{equation}\label{6.21}
  \|u\|_{\alpha,I}\le\|u\|_{0,I}+(b-a)^{1-\alpha}\|u\|_{1,I}\,.
\end{equation}
\end{proof}
\begin{lemma}
(Estimate of the $C^{1+\alpha}$-norm)\\[0.1cm]
Let $u\in C^{2+\alpha}(I,\mathbb R)$ solve (6.1b). Then, for all real $0<\mu<\frac{1}{2}$ it holds
\begin{equation}\label{6.22}
  \|u\|_{1+\alpha,I}\le\|u\|_{1,I}+c_1\|u\|_{0,I}+c_2\|f\|_{0,I}+\mu c_3\|u\|_{2+\alpha,I}
\end{equation}
with the above bounds on $\|u\|_{1,I}$ and the constants
\begin{equation}\label{6.23}
\begin{array}{l}
  \displaystyle
  c_1
  :=(b-a)^{1-\alpha}
        \left(
          \frac{p_1'}{p_0}
          +\frac{2p_1p_1'}{\mu p_0^2(b-a)}
          +\frac{p_1'q_1(b-a)}{2p_0^2}
          +\frac{q_1}{p_0}
        \right), \\[0.7cm]
  \displaystyle
  c_2
  :=(b-a)^{1-\alpha}
        \left(
          \frac{p_1'(b-a)}{2p_0^2}+\frac{1}{p_0}
        \right),\quad
  c_3:=\frac{p_1p_1'(b-a)^{2-\alpha}}{2p_0^2}\,.
\end{array}
\end{equation}
\end{lemma}
\begin{proof}
Let $x_1,x_2\in I.$ We find $\xi\in[x_1,x_2]$ with the property
\begin{equation}\label{6.24}
  \frac{|u'(x_1)-u'(x_2)|}{|x_1-x_2|^\alpha}
  =|u''(\xi)||x_1-x_2|^{1-\alpha}
  \le (b-a)^{1-\alpha}[u]_{2,I}\,.
\end{equation}
From the equation
\begin{equation}\label{6.25}
  u''+\frac{p'}{p}\,u'+\frac{q}{p}\,u=\frac{1}{p}\,f
\end{equation}
we conclude
\begin{equation}\label{6.26}
  [u]_{2,I}\le\frac{p_1'}{p_0}\,\|u\|_{1,I}+\frac{q_1}{p_0}\,\|u\|_{0,I}+\frac{1}{p_0}\,\|f\|_{0,I}\,.
\end{equation}
Using the above $C^1$-estimate we get
\begin{equation}\label{6.27}
  [u]_{2,I}\le c_1\|u\|_{0,I}+c_2\|f\|_{0,I}+c_3\|u\|_{2+\alpha,I}
\end{equation}
with the given constants $c_1,c_2,c_3.$
\end{proof}
\noindent
Let us define
\begin{equation}\label{6.28}
  c_4:=c_2+\frac{b-a}{2p_0}\,,\quad
  c_5:=1+c_1+\frac{2p_1}{\mu p_0(b-a)}+\frac{q_1(b-a)}{2p_0}\,,\quad
  c_6:=c_3+\frac{p_1(b-a)}{2p_0}\,.
\end{equation}
From (\ref{6.22}) it follows that
\begin{equation}\label{6.29}
  \|u\|_{1+\alpha,I}
  \le c_4\|f\|_{0,I}+c_5\|u\|_{0,I}+\mu c_6\|u\|_{2+\alpha,I}\,.
\end{equation}
To prove global estimates of the $C^{2+\alpha}$-norm we note that
\begin{equation}\label{6.30}
  u''(x)=\frac{f(x)}{p(x)}-\frac{p'(x)}{p(x)}\,u'(x)-\frac{q(x)}{p(x)}\,u(x)\equiv F(x),
\end{equation}
that is, we have to investigate the Poisson-type equation $u''=F.$
\begin{lemma}
(Estimate of the $C^0$-norm for Poisson's equation)\\[0.1cm]
Let $u\in C^{2+\alpha}(I,\mathbb R)$ solve (\ref{6.30}). Then
\begin{equation}\label{6.31}
  \|u\|_{0,I}\le 3\max_{x\in\partial I}|u(x)|+\|F\|_{0,I}e^{b-2a}\,.
\end{equation}
\end{lemma}
\noindent
This follows from the $C^0$-estimate (\ref{6.7}) with $p\equiv 1,$ $q\equiv 0.$
\begin{lemma}
(Estimate of the $C^1$-norm for Poisson's equation)\\[0.1cm]
Let $u\in C^{2+\alpha}(I,\mathbb R)$ solve (\ref{6.30}). For all real $0<\mu<\frac{1}{2}$ it holds
\begin{equation}\label{6.32}
  \|u\|_{1,I}
  \le\frac{(1+\mu)(b-a)}{2}\,\|F\|_{0,I}+\frac{2+\mu(b-a)}{\mu(b-a)}\,\|u\|_{0,I}\,.
\end{equation}
\end{lemma}
\begin{proof}
As in (\ref{6.15}) we have
\begin{equation}\label{6.33}
  |u'(c)|
  \le\frac{2}{\mu(b-a)}\,\|u\|_{0,I}+\frac{\mu(b-a)}{2}\,[u]_{2,I}
  =\frac{2}{\mu(b-a)}\,\|u\|_{0,I}+\frac{\mu(b-a)}{2}\,\|F\|_{0,I}
\end{equation}
because $[u]_{2,I}=\|F\|_{0,I}.$ Integrating $u''(x)=F(x)$ from $\xi=c$ to $\xi=x$ gives
\begin{equation}\label{6.34}
  [u]_{1,I}
  \le\frac{b-a}{2}\,\|F\|_{0,I}+\frac{2}{\mu(b-a)}\,\|u\|_{0,I}+\frac{\mu(b-a)}{2}\,\|F\|_{0,I}\,.
\end{equation}
This proves the statement.
\end{proof}
\noindent
Furthermore, $[u'']_{2,I}=\|F\|_{0,I}$ yields immediately
\begin{lemma}
(Estimate of the $C^2$-norm for Poisson's equation)\\[0.1cm]
Let $u\in C^{2+\alpha}(I,\mathbb R)$ solve (\ref{6.30}). For all real $0<\mu<\frac{1}{2}$ it holds
\begin{equation}\label{6.35}
  \|u\|_{2,I}
  \le\frac{2+(1+\mu)(b-a)}{2}\,\|F\|_{0,I}+\frac{2+\mu(b-a)}{\mu(b-a)}\,\|u\|_{0,I}\,.
\end{equation}
\end{lemma}
\noindent
As in (\ref{6.18}) and (\ref{6.22}) we prove
\begin{lemma}
(Estimate of the $C^\alpha$-norm for Poisson's equation)\\[0.1cm]
Let $u\in C^{2+\alpha}(I,\mathbb R)$ solve (\ref{6.30}). Then it holds
\begin{equation}\label{6.36}
  \|u\|_{\alpha,I}\le\|u\|_{0,I}+(b-a)^{1-\alpha}[u]_{1,I}
\end{equation}
with the above estimates of $\|u\|_{0,I}$ and $[u]_{1,I}.$
\end{lemma}
\begin{lemma}
(Estimate of the $C^{1+\alpha}$-norm for Poisson's equation)\\[0.1cm]
Let $u\in C^{2+\alpha}(I,\mathbb R)$ solve (\ref{6.30}). Then it holds
\begin{equation}\label{6.37}
  \|u\|_{1+\alpha,I}
  \le\|u\|_{1,I}+(b-a)^{1-\alpha}[u]_{2,I}
\end{equation}
with the above bounds on $\|u\|_{1,I}$ and $[u]_{2,I}.$
\end{lemma}
\begin{lemma}
(Estimate of the $C^{2+\alpha}$-norm for Poisson's equation)\\[0.1cm]
Let $u\in C^{2+\alpha}(I,\mathbb R)$ solve (\ref{6.30}). Then it holds
\begin{equation}\label{6.38}
  \|u\|_{2+\alpha,I}=\|u\|_{2,I}+[u]_{2+\alpha,I}=\|u\|_{2,I}+\|F\|_{\alpha,I}
\end{equation}
with the above bound on $\|u\|_{2,I}.$
\end{lemma}
\noindent
Our main result is a global estimate of the $C^{2+\alpha}$-norm for regular solutions of (6.1b).
\begin{proposition}
(Global estimate of the $C^{2+\alpha}$-norm)\\[0.1cm]
Let $u\in C^{2+\alpha}(I,\mathbb R)$ solve (6.1b). Let
\begin{equation}\label{6.39}
  \widehat p_1:=\|p\|_{\alpha,I}\,,\quad
  \widehat q_1:=\|q\|_{\alpha,I}\,,\quad
  \widehat p_1\!':=\|p\|_{1+\alpha,I}\,.
\end{equation}
Then there exist real constants $C_1,C_2=C_1,C_2(p_0,\widehat p_1,\widehat q_1,\widehat p_1\!',b-a,\alpha)\in(0,+\infty)$ such that
\begin{equation}\label{6.40}
  \|u\|_{2+\alpha,I}
  \le C_1\|f\|_{\alpha,I}+C_2\|u\|_{0,I}\,.
\end{equation}
\end{proposition}
\begin{proof}
Note that
\begin{equation}\label{6.41}
  u''(x)=-\frac{p'(x)}{p(x)}\,u'(x)-\frac{q(x)}{p(x)}\,u(x)+\frac{f(x)}{p(x)}\equiv F(x).
\end{equation}
From (\ref{6.35}) it follows that
\begin{equation}\label{6.42}
\begin{array}{lll}
  \|u\|_{2+\alpha,I}\negthickspace\displaystyle
  &  =  & \negthickspace\displaystyle
          \|u\|_{2,I}+[u]_{2+\alpha,I} 
          \,\le\,\|u\|_{2,I}+\|F\|_{\alpha,I} \\[0.4cm]
  & \le & \negthickspace\displaystyle
          \frac{2+(1+\mu)(b-a)}{2}\,\|F\|_{0,I}+\|F\|_{\alpha,I}+\frac{2+\mu(b-a)}{\mu(b-a)}\,\|u\|_{0,I}\,.
\end{array}
\end{equation}
\begin{itemize}
\item[1.]
We estimate $\|F\|_{0,I}.$ From (\ref{6.29}) we infer
\begin{equation}\label{6.43}
\begin{array}{lll}
  \|F\|_{0,I}\negthickspace
  & \le & \negthickspace\displaystyle
          \frac{1}{p_0}\,\|f\|_{0,I}+\frac{p_1'}{p_0}\,\|u\|_{1,I}+\frac{q_1}{p_0}\,\|u\|_{0,I} \\[0.6cm]
  & \le & \negthickspace\displaystyle
          \frac{1}{p_0}\,(1+p_1'c_4)\|f\|_{0,I}
          +\frac{1}{p_0}\,(p_1'c_5+q_1)\|u\|_{0,I}
          +\mu\frac{p_1'}{p_0}\,c_6\|u\|_{2+\alpha,I}\,.
\end{array}
\end{equation}
We set
\begin{equation}\label{6.44}
  c_7:=\frac{1}{p_0}\,(1+p_1'c_4),\quad
  c_8:=\frac{1}{p_0}\,(p_1'c_5+q_1),\quad
  c_9:=\frac{p_1'}{p_0}\,c_6\,,
\end{equation}
such that
\begin{equation}\label{6.45}
  \|F\|_{0,I}\le c_7\|f\|_{0,I}+c_8\|u\|_{0,I}+\mu c_9\|u\|_{2+\alpha,I}\,.
\end{equation}
\item[2.]
Analogously we estimate
\begin{equation}\label{6.46}
  \|F\|_{\alpha,I}
  \le\frac{\|f\|_{\alpha,I}}{p_0}
     +\frac{\widehat p_1\!'}{p_0}\,\|u\|_{1+\alpha,I}
     +\frac{\widehat q_1}{p_0}\,\|u\|_{\alpha,I}\,.
\end{equation}
We define
\begin{equation}\label{6.47}
\begin{array}{lll}
  c_{10}\negthickspace
  & := & \negthickspace\displaystyle
         \frac{1}{p_0}\,\big(1+\widehat p_1\!'c_4+\widehat q_1(b-a)^{1-\alpha}c_4\big), \\[0.6cm]
  c_{11}\negthickspace
  & := & \negthickspace\displaystyle
         \frac{1}{p_0}\,\big(\widehat p_1\!'c_5+\widehat q_1+\widehat q_1(b-a)^{1-\alpha}c_5\big), \\[0.6cm]
  c_{12}\negthickspace
  & := & \negthickspace\displaystyle
         \frac{1}{p_0}\big(\widehat p_1\!'c_6+\widehat q_1(b-a)^{1-\alpha}c_6\big).
\end{array}
\end{equation}
Together with (\ref{6.18}) and (\ref{6.22}) we have
\begin{equation}\label{6.48}
  \|F\|_{\alpha,I}\le c_{10}\|f\|_{\alpha,I}+c_{11}\|u\|_{0,I}+\mu c_{12}\|u\|_{2+\alpha,I}\,.
\end{equation}
\end{itemize}
Now, we set
\begin{equation}\label{6.49}
\begin{array}{lll}
  c_{13}\negthickspace
  & := & \negthickspace\displaystyle
         \frac{2+(1+\mu)(b-a)}{2}\,c_7+c_{10}\,, \\[0.6cm]
  c_{14}\negthickspace
  & := & \negthickspace\displaystyle
         \frac{2+(1+\mu)(b-a)}{2}\,c_8+\frac{2+\mu(b-a)}{\mu(b-a)}+c_{11}\,, \\[0.6cm]
  c_{15}\negthickspace
  & := & \negthickspace\displaystyle
         \frac{2+(1+\mu)(b-a)}{2}\,c_9+c_{12}\,.
\end{array}
\end{equation}
Summarizing (\ref{6.42}), (\ref{6.45}) and (\ref{6.48}) we arrive at the estimate
\begin{equation}\label{6.50}
  \|u\|_{2+\alpha,I}\le c_{13}\|f\|_{\alpha,I}+c_{14}\|u\|_{0,I}+\mu c_{15}\|u\|_{2+\alpha,I}\,.
\end{equation}
Finally, we choose $0<\mu<\frac{1}{2}$ such that $1-\mu c_{15}>0.$ It follows that
\begin{equation}\label{6.51}
  \|u\|_{2+\alpha,I}
  \le\frac{c_{13}}{1-\mu c_{15}}\,\|f\|_{\alpha,I}
     +\frac{c_{14}}{1-\mu c_{15}}\,\|u\|_{0,I}\,.
\end{equation}
Setting $C_1:=c_{13}(1-\mu c_{15})^{-1}$ and $C_2:=c_{14}(1-\mu c_{15})^{-1}$ proves the statement.
\end{proof}
\noindent
The presented a priori estimates are not sharp. In the foregoing application we needed only a qualitative form of these constants.

\vspace*{1.6cm}
Steffen Fr\"ohlich\\
Technische Universit\"at Darmstadt\\
Fachbereich Mathematik, AG 4\\
Differentialgeometrie und Geometrische Datenverarbeitung\\
Schlo\ss{}gartenstra\ss{}e 7\\
D-64289 Darmstadt\\
Germany\\[0.1cm]
e-mail: sfroehlich@mathematik.tu-darmstadt.de

\end{document}